# SPECTRAL ESTIMATION OF THE FRACTIONAL ORDER OF A LÉVY PROCESS

By Denis Belomestny[1]

*Weierstrass-Institute Berlin*

We consider the problem of estimating the fractional order of a Lévy process from low frequency historical and options data. An estimation methodology is developed which allows us to treat both estimation and calibration problems in a unified way. The corresponding procedure consists of two steps: the estimation of a conditional characteristic function and the weighted least squares estimation of the fractional order in spectral domain. While the second step is identical for both calibration and estimation, the first one depends on the problem at hand. Minimax rates of convergence for the fractional order estimate are derived, the asymptotic normality is proved and a data-driven algorithm based on aggregation is proposed. The performance of the estimator in both estimation and calibration setups is illustrated by a simulation study.

**1. Introduction.** Nowadays Lévy processes are undoubtedly one of the most popular tool for modeling economic and financial time series [see, e.g., Cont and Tankov (2004), for an overview]. This is not surprising if one takes into account their simplicity and analytic tractability on the one hand and the ability to reproduce many stylized facts of financial time series on the other hand. In the last decade, new subclasses of Lévy processes have been introduced and actively studied (mainly in the context of option pricing). Among the best known models are normal inverse Gaussian processes (NIG), hyperbolic processes (HP), generalized hyperbolic processes (GHP) and truncated (or tempered) Lévy processes (TLP). Boyarchenko and Levendorskiĭ (2002) have introduced a general class of regular Lévy processes of exponential type (RLE) which contains all above mentioned particular Lévy models. This type of processes is characterized

Received August 2008; revised May 2009.
[1]Supported in part by SFB 649 "Economic Risk."

by the requirement that the modulus of the characteristic function of increments behaves like $\exp(-\eta|u|^\alpha)$ as $|u| \to \infty$ for some $0 < \alpha < 2$. Parameter $\alpha$ coincides with the fractional order of the underlying Lévy process and plays an important role because it determines the decay of the characteristic function and hence the smoothness properties of the corresponding state price density. Statistical inference for RLE processes is the subject of our paper.

There are basically two types of statistical problems relevant for Lévy processes: the estimation of parameters of a Lévy process $X_t$ from a time series of the asset $S_t = \exp(X_t)$ and the calibration of these parameters using options data. Both problems have received much attention recently.

Suppose that a Lévy process $X_t$ is observed at $n$ time points $\Delta, 2\Delta, \ldots, n\Delta$. Since $X_0 = 0$, this amounts to observing $n$ increments $\chi_i = X_{i\Delta} - X_{(i-1)\Delta}$, $i = 1, \ldots, n$. If $\Delta$ is small (high-frequency data), then a large increment $\chi_i$ indicates that a jump occurred between time $t_{i-1}$ and $t_i$. Based on this insight and the continuous-time observation analogue, inference for the Lévy measure of the underlying Lévy process can be conducted. See, for example, Aït-Sahalia and Jacod (2006) for a semiparametric problem of estimating volatility of a stable process under the presence of Lévy perturbation or Lee and Mykland (2008) and Figueroa-López and Houdré (2006) for the nonparametric problem of testing and estimation for jump diffusion models. For low-frequency observations, however, we cannot be sure to what extent the increment $\chi_i$ is due to one or several jumps or just to the diffusion part of the Lévy process. The only way to draw inference is to use the fact that the increments form independent realizations of infinitely divisible probability distributions. In this setting, a variety of methods have been proposed in the literature: standard maximum likelihood estimation DuMouchel (1973a, 1973b, 1975), using the empirical characteristic function as an estimating equation [see, e.g., Press (1972), Fenech (1976), Feuerverger and McDunnough (1981a), Singleton (2001)], maximum likelihood by Fourier inversion of the characteristic function Feuerverger and McDunnough (1981b), a regression based on the explicit form of the characteristic function Koutrouvelis (1980) or other numerical approximations Nolan (1997). Some of these methods were compared in Akgiray and Lamoureux (1989). Note that all of the aforementioned papers deal with the specific parametric (mainly stable) models. A semiparametric estimation problem for Lévy models has recently been considered in Neumann and Reiss (2009) and Gugushvili (2008).

The second calibration problem is of special importance for financial applications because pricing of options is performed under an equivalent martingale measure, and one can infer on this measure only from options data. Since option data is sparse and the underlying inverse problem is usually ill-posed, we face a rather complicated estimation issue. Different approaches have been proposed in the literature to regularize the underlying inverse problem. For example, in Cont and Tankov (2004) and Cont and Tankov (2006), a method based on the penalized least squares estimation with the minimal entropy penalization is proposed. Belomestny and Reiss (2006) developed a spectral calibration method which avoids solving a high-dimensional optimization problem and is based on the direct inversion of a Fourier pricing formula with a cut-off regularization in spectral domain. This method essentially employees the integrability property of the underlying Lévy measure (finite activity Lévy processes) that excludes many interesting infinite activity Lévy processes.

In this paper we consider the problem of estimating the fractional order of a Lévy process from low-frequency historical as well as options data. Our problem is a semiparametric one because we do not assume any specific parametric model for the underlying process but only some asymptotic behavior. The spectral approach allows us to treat both estimation and calibration problems in a unified framework and leads to an efficient data-driven algorithm. Moreover, the fractional order estimate delivered by the spectral method possesses several interesting optimality properties.

The problem of estimating the degree of activity of jumps in semimartingale framework using high-frequency financial data has recently been considered in Aït-Sahalia and Jacod (2009). On the one hand, small increments of the process turn out to be most informative for estimating the activity index. On the other hand, these small increments are the ones where the contribution from the continuous martingale part is mixed with the contribution from the small jumps. Aït-Sahalia and Jacod (2009) proposed an estimation procedure which is able to "see through" the continuous part and consistently estimate the degree of activity for the small jumps under some restrictions on the structure of the underlying semimartingale. Note that in the case of Lévy processes the degree of activity of jumps is identical to the fractional order of the underlying Lévy process. We also stress that the case when both diffusion and jump components are presented can be treated in the framework of spectral estimation as well (see Section 6.9).

Short outline of the paper. In Section 2 we introduce the class of RLE processes. Section 3 discusses some aspects of financial modeling with RLE processes. Section 4 describes the observational model. In Section 5 methods of estimating the characteristic function of a Lévy process from low-frequency historical and options data are presented. Section 6 is devoted to the spectral calibration method of estimating the fractional order $\alpha$. We discuss here the problems of regularization and derive minimax rates of convergence for a class of Lévy processes. In Section 7 an adaptive procedure for estimating $\alpha$ is presented, and its properties are discussed. We conclude with some simulation results.

## 2. Regular Lévy processes of exponential type.

In this section we recall some basic properties of Lévy processes.

2.1. *Spectral properties of Lévy processes.* Consider a Lévy process $X_t$ with a Lévy measure $\nu$. That is, $X_t$ is càdlàg process with independent and stationary increments such that the characteristic function of its marginals $\phi_t(u)$ is given by

$$\begin{aligned}\phi_t(u) &:= \mathrm{E}[e^{\mathrm{i}uX_t}] \\ &= \exp\left\{t\left(\mathrm{i}u\mu - \frac{u^2a^2}{2} + \int_{\mathbb{R}}(e^{\mathrm{i}ux} - 1 - \mathrm{i}ux1_{\{|x|\le 1\}})\nu(dx)\right)\right\}.\end{aligned} \tag{2.1}$$

So, any Lévy process $X_t$ is characterized by the so called Lévy triple $(\mu, a, \nu)$ where $\mu \in \mathbb{R}$ is a drift, $a > 0$ is a diffusion volatility and $\nu$ is a Lévy measure. Note that the drift $\mu$ depends on the type of truncation in (2.1). In fact, this characterization is unique for a fixed truncation function and we can reconstruct the Lévy triple from the characteristic function $\phi_t(u)$. This reconstruction may be viewed as consisting of three steps. First, because of

$$\frac{1}{|u|^2}\int_{\mathbb{R}}(e^{\mathrm{i}ux} - 1 - \mathrm{i}ux1_{\{|x|\le 1\}})\nu(dx) \to 0, \qquad |u| \to \infty, \tag{2.2}$$

we can find $a^2/2$ as $\lim_{|u|\to\infty}|u|^{-2}\psi(u)$ with

$$\psi(u) = t^{-1}\log(\phi_t(u)).$$

Second, note that

$$\int_{-1}^{1}(\widetilde{\psi}(u) - \widetilde{\psi}(u+w))\,dw = \int_{\mathbb{R}} e^{\mathrm{i}ux}\rho(dx)$$

with

$$\widetilde{\psi}(u) = \psi(u) + \frac{a^2}{2}u^2, \qquad \rho(dx) = 2\left(1 - \frac{\sin x}{x}\right)\nu(dx).$$

Since $\rho$ is a finite measure ($\int(x^2 \wedge 1)\nu(dx) < \infty$), one can uniquely reconstruct it (and hence $\nu$) from $\widetilde{\psi}(u)$. Finally, we find $\mu$ as $\lim_{u\to\infty}[\widetilde{\psi}(u)/(\mathrm{i}u)]$. So, in principle, we can recover all characteristics of the underlying Lévy process (including the fractional order) provided that $\phi_t$ is completely known. If, however, $\phi_t$ is estimated from data we face an ill-posed estimation problem because a small perturbation in $\phi_t$ may deteriorate its asymptotic behavior and lead to the violation of (2.2). In this case using a regularization technique [see, e.g., Cont and Tankov (2004) or Belomestny and Reiss (2006)], we still can get an asymptotically consistent estimates for the whole triple $(\mu, a, \nu)$ given a consistent estimate of $\phi_t$.

REMARK 2.1. A consistent estimation of $\psi(u)$ from a time series of $X_t$ is only possible if the number of observations from the distribution with the cf. $\phi_t(u)$ for some $t>0$ increases. This can be either due to a decreasing time step in a times series of the process $X$ (high frequency data) or due to an increasing time horizon (low frequency data). While the first type of observational models has received much attention in recent years, there are only few papers dealing with low frequency data [see, e.g., Neumann and Reiss (2009)].

2.2. *Fractional order of Lévy processes.* Let $X_t$ be a Lévy process with a Lévy measure $\nu$. The value

$$\alpha := \inf\left\{r \geq 0 : \int_{|x|\leq 1} |x|^r \nu(dx) < \infty\right\}$$

is called the *fractional order* or the Blumenthal–Getoor index of the Lévy process $X_t$. This index $\alpha$ is related to the "degree of activity" of jumps. All Lévy measures put finite mass on the set $(-\infty,-\epsilon]\cup[\epsilon,\infty)$ for any arbitrary $\epsilon>0$, so if the process has infinite jump activity, it must be because of the small "jumps," defined as those smaller than $\epsilon$. If $\nu([-\epsilon,\epsilon])<\infty$ the process has finite activity and $\alpha=0$. But if $\nu([-\epsilon,\epsilon])=\infty$, that is, the process has infinite activity and in addition the Lévy measure $\nu((-\infty,-\epsilon]\cup[\epsilon,\infty))$ diverges near 0 at a rate $|\epsilon|^{-\alpha}$ for some $\alpha>0$, then the fractional order of $X_t$ is equal to $\alpha$. The higher $\alpha$ gets, the more frequent the small jumps become [see Aït-Sahalia and Jacod (2009) for more discussion].

The Blumenthal–Getoor index is closely related to the notion of the degree of jump activity that applies to general semimartingales as shown in Aït-Sahalia and Jacod (2009), and reduces to the Blumenthal–Getoor index in the special case of Lévy processes.

Note also that the Blumenthal–Getoor index coincides with the stability index for stable processes. Another example of processes having a prescribed fractional order $\alpha$ is the class of tempered stable processes of order $\alpha$. Boyarchenko and Levendorskiĭ (2002) studied a generalization of tempered stable processes, called *regular Lévy processes of exponential type* (RLE). A Lévy process is said to be a RLE process of type $[\lambda_-,\lambda_+]$ and order $\alpha\in(0,2)$ if the Lévy measure has exponentially decaying tails with rates $\lambda_-\geq 0$ and $\lambda_+\geq 0$

$$\int_{-\infty}^{-1} e^{\lambda_-|y|}\nu(dy)<\infty, \qquad \int_1^{\infty} e^{\lambda_+ y}\nu(dy)<\infty \tag{2.3}$$

and behaves near zero as $|y|^{-(1+\alpha)}$;

$$\int_{|y|>\epsilon}\nu(dy) \asymp \frac{\Pi(\epsilon)}{\epsilon^{\alpha}}, \qquad \epsilon\to+0,$$

where $\Pi$ is some positive function on $\mathbb{R}_+$ satisfying $0<\Pi(+0)<\infty$. Obviously, the fractional order of an RLE process of order $\alpha$ is equal to $\alpha$. An equivalent definition of an RLE process in terms of its characteristic exponent $\psi(u)$ can be given as follows. A Lévy process is considered to be an RLE process of type $[\lambda_-, \lambda_+]$ and order $\alpha \in (0,2)$ if the following representation holds:

$$\psi(u) = \mathrm{i}\mu u + \vartheta(u), \qquad \mu \in \mathbb{R}, \tag{2.4}$$

where function $\vartheta$ admits a continuation from $\mathbb{R}$ into the strip $\{z \in \mathbb{C} : \operatorname{Im} z \in [-\lambda_+, \lambda_-]\}$ and is of the form

$$\vartheta(u) = -|u|^{\alpha}\pi(u), \tag{2.5}$$

where $\pi(u)$ is a function satisfying

$$\limsup_{|u|\to\infty} |\pi(u)| < \infty \quad \text{and} \quad \liminf_{|u|\to\infty} |\pi(u)| > 0$$

such that

$$\operatorname{Re}[\pi(u)] > 0, \qquad u \in \mathbb{R} \setminus \{0\}. \tag{2.6}$$

As was mentioned in the Introduction, the class of RLE processes includes among others hyperbolic, normal inverse Gaussian and tempered stable processes but does not include variance Gamma process. In the sequel we will mainly consider RLE processes without regularity conditions (2.3) (or equivalently with $\lambda_- = \lambda_+ = 0$) since only the behavior of a Lévy measure near zero matters for the fractional order of the corresponding Lévy process.

As mentioned before, in this work we are going to consider the problem of estimating the fractional order $\alpha$ of a Lévy process $X_t$ from a time series of asset prices as well as from option prices. Before turning to this, let us first make our modelling and observational framework more precise.

**3. Financial modelling.** In this section we recall basic facts concerning financial modelling with exponential Lévy models.

3.1. *Asset dynamics.* We assume that the asset price $S_t$ follows an exponential Lévy model under both *historical measure* $\mathbb{P}$ and *risk neutral measure* $\mathbb{Q}$. Specifically, we suppose that

$$S_t = \begin{cases} S e^{X_t}, & \text{under } \mathbb{P}, \\ S e^{rt + Y_t}, & \text{under } \mathbb{Q}, \end{cases}$$

where $X_t$ and $Y_t$ are Lévy processes, $S > 0$ is the present value of the asset (at time 0) and $r \geq 0$ is the riskless interest rate which is assumed to be known and constant. Note that the martingale condition for $S_t$ under $\mathbb{Q}$ entails $\mathrm{E}^{\mathbb{Q}}[e^{Y_t}] = 1$. The martingale measure $\mathbb{Q}$ is in fact not unique under the presence of jumps. As is standard in the calibration literature, it is assumed to be settled by the market and to be identical for all options under consideration. Processes $X_t$ and $Y_t$ are related by the requirement that measures $\mathbb{P}$ and $\mathbb{Q}$ ought to be equivalent: $\mathbb{P} \sim \mathbb{Q}$. Interestingly, this implies that if $X_t$ and $Y_t$ are RLE process and $X_t$ is of order $\alpha^{\mathbb{P}}$, then $Y_t$ has the order $\alpha^{\mathbb{Q}} = \alpha^{\mathbb{P}}$. Indeed, the equivalence of the corresponding Lévy measures $\nu^{\mathbb{P}}$ and $\nu^{\mathbb{Q}}$ implies [see Sato (1999)]

$$\int_0^\infty (\sqrt{d\nu^{\mathbb{Q}}/d\nu^{\mathbb{P}}} - 1)^2 \nu^{\mathbb{P}}(dx) < \infty. \tag{3.1}$$

Since for RLE processes $d\nu^{\mathbb{Q}}(x)/d\nu^{\mathbb{P}}(x) \asymp x^{(\alpha^{\mathbb{P}} - \alpha^{\mathbb{Q}})}$ and $d\nu^{\mathbb{P}}(x) \asymp x^{-(1+\alpha^{\mathbb{P}})}\,dx$ as $x \to +0$, the condition (3.1) can be satisfied only if $\alpha^{\mathbb{P}} = \alpha^{\mathbb{Q}}$. This means that the fractional order of the underlying Lévy process must be the same under both historical and risk-neutral measures. This not only indicates the importance of the fractional order parameter for financial applications but also suggests that the combination of two estimates of the fractional order $\alpha$ under $\mathbb{P}$ and $\mathbb{Q}$ might be useful, for example, to reduce the overall variance of the resulting combined estimator.

3.2. *Option pricing.* The risk neutral price at time $t = 0$ of the European call option with strike $K$ and maturity $T$ is given by

$$C(K, T) = e^{-rT} \mathrm{E}^{\mathbb{Q}}[(S_T - K)^+].$$

Using the independence of increments, we can reduce the number of parameters by introducing the so-called negative log-forward moneyness,

$$y := \log(K/S) - rT,$$

such that the call price in terms of $y$ is given by

$$\mathcal{C}(y, T) = S\mathrm{E}^{\mathbb{Q}}[(e^{Y_T} - e^y)^+].$$

The analogous formula for the price of the European put option is $\mathcal{P}(y, T) = S\mathrm{E}^{\mathbb{Q}}[(e^y - e^{Y_T})^+]$, and a well-known put-call parity is easily established;

$$\mathcal{C}(y, T) - \mathcal{P}(y, T) = S\mathrm{E}^{\mathbb{Q}}[e^{Y_T} - e^y] = S(1 - e^y).$$

As we need to employ Fourier techniques, we introduce the function

$$O_T(y) := \begin{cases} S^{-1}\mathcal{C}(y, T), & y \geq 0, \\ S^{-1}\mathcal{P}(y, T), & y < 0. \end{cases} \tag{3.2}$$

The function $O_T$ records normalized call prices for $y \geq 0$ and normalized put prices for $y < 0$. It possesses many interesting properties [see Belomestny and Reiss (2006) for details]; one of them being the following connection between the Fourier transform of $O_T$ and the characteristic function of $Y_T$ denoted by $\phi_T^{\mathbb{Q}}$:

$$\mathbf{F}[O_T](v) = \frac{1 - \phi_T^{\mathbb{Q}}(v - \mathrm{i})}{v(v - \mathrm{i})}, \qquad v \in \mathbb{R}. \tag{3.3}$$

Another property which directly follows from (3.3) is that

$$\int_{\mathbb{R}} e^{-2y} O_T(y)\,dy < \infty, \tag{3.4}$$

provided that $\mathrm{E}[e^{2Y_T}]$ exists and is finite.

**4. Observations.** We consider two kinds of observational models corresponding to two types of statistical problems we are going to tackle. While the first type of models assumes the time series of $S_t$ is directly available, the second one supposes that only some functionals of $S_t$ can be observed.

4.1. *Time series data.* We assume that the values of the log-price process $X_t = \log(S_t)$ on equidistant time grid $\pi = \{t_0, t_1, \ldots, t_n\}$ are observed.

4.2. *Option data.* As to option data, we assume that we will be given the prices of $n$ call options for a set of forward log-moneynesses $y_0 < y_1 < \cdots < y_n$ and a fixed maturity $T$ corrupted by noise. In terms of the function $O$, the following sample is available:

$$\mathcal{O}_T(y_j) = O_T(y_j) + \sigma(y_j)\xi_j, \qquad j = 1, \ldots, n. \tag{4.1}$$

It is supposed that $\{\xi_j\}$ are independent, centered, random variables with $\mathrm{E}[\xi_j^2] = 1$ and $\sup_j \mathrm{E}[\xi_j^4] < \infty$. Furthermore, we assume that

$$\int_{\mathbb{R}} e^{-2y} \sigma^2(y)\,dy < \infty.$$

This condition is required because we need to transform the original regression model (4.1) to an exponentially weighted one,

$$\widetilde{\mathcal{O}}_T(y_j) = \widetilde{O}_T(y_j) + \widetilde{\sigma}(y_j)\xi_j, \qquad j = 1, \ldots, n, \tag{4.2}$$

with $\widetilde{\mathcal{O}}_T(y) = e^{-y}\mathcal{O}_T(y)$, $\widetilde{O}_T(y) = e^{-y}O_T(y)$ and $\widetilde{\sigma}(y) = e^{-y}\sigma(y)$.

As a matter of fact, a consistent estimation of the fractional order $\alpha$ is only possible if the amount of data available increases. In our asymptotic analysis we will therefore assume that the number of time series observations and the number of available options tend to infinity.

**5. Estimation of characteristic functions $\phi^{\mathbb{P}}$ and $\phi^{\mathbb{Q}}$.** The main idea of the spectral estimation method (SEM) is to infer on the parameters of the underlying model using its special structure in the spectral domain. Since spectral behavior of a RLE process is described explicitly by (2.4) and (2.5), we can apply SEM as soon as an estimate for the corresponding characteristic function is available. While estimation of $\phi$ under $\mathbb{P}$ is rather straightforward, its calibration from option prices under $\mathbb{Q}$ requires special treatment.

5.1. *Estimation of $\phi$ under $\mathbb{P}$.* We estimate the characteristic function $\phi^{\mathbb{P}}_{|\pi|}(u)$ by its empirical counterpart,

$$\widetilde{\phi}^{\mathbb{P}}_{|\pi|}(u) = \frac{1}{n}\sum_{j=1}^{n} e^{\mathrm{i}u(X_{t_j}-X_{t_{j-1}})}.$$

The *empirical characteristic function* $\widetilde{\phi}^{\mathbb{P}}_{|\pi|}$ possesses many interesting properties, and we refer to Ushakov (1999) for a comprehensive overview.

5.2. *Estimation of $\phi$ under $\mathbb{Q}$.* For estimating $\phi^{\mathbb{Q}}_T$ we employ the Fourier technique. So, motivated by (3.3), we define

$$\widetilde{\phi}^{\mathbb{Q}}_T(u) := 1 - u(u+\mathrm{i})\left[\sum_{j=1}^{n}\delta_j \widetilde{\mathcal{O}}_T(y_j)e^{\mathrm{i}uy_j}\right], \qquad u\in\mathbb{R}, \tag{5.1}$$

where $\delta_j = y_j - y_{j-1}$ and $\widetilde{\mathcal{O}}_T$ is defined in (4.2). For more involved methods of approximating $\mathbf{F}[O_T](u)$ see Belomestny and Reiss (2006).

**6. Estimation of fractional order.** In this section we turn to the problem of estimating the fractional order of a RLE process. To this aim we apply the spectral estimation method accompanied by a spectral cut-off regularization.

6.1. *Main idea.* Let us consider a RLE process with the characteristic exponent $\psi(u)$ of the form (2.4) and (2.5). In the sequel we assume (mainly for the sake of simplicity) that $\lim_{u\to-\infty}\pi(u) = \lim_{u\to\infty}\pi(u) = \eta\in\mathbb{R}_+$. In this case we can rewrite $\vartheta$ as

$$\vartheta(u) = -\eta|u|^{\alpha}\tau(u), \tag{6.1}$$

where $\mathrm{Re}[\tau(u)] > 0$ for $u\in\mathbb{R}\setminus\{0\}$ and $\tau(u)\to 1$ as $|u|\to\infty$. The formula,

$$\begin{aligned}\mathcal{Y}(u) &:= \log(-\log(|\phi(u)|^2)) \\ &= \log(2\eta) + \alpha\log(u) + \log(\mathrm{Re}\,\tau(u)), \qquad u>0,\end{aligned} \tag{6.2}$$

with $\phi(u) = \exp(\psi(u))$, suggests how to estimate $\alpha$ from $\phi$. Indeed, in terms of the new "data" $\mathcal{Y}$, we have a linear semiparametric problem with the "nuisance" nonparametric part $\log(\mathrm{Re}\,\tau(u))$. Since $\log(\mathrm{Re}\,\tau(u))$ tends to 0 as $|u|\to\infty$, we can get rid of this component by basing our estimation on $\mathcal{Y}(u)$ with large $|u|$. On the other hand, if we plug-in an estimate $\widetilde{\phi}$ instead of $\phi$, the variance of $\mathcal{Y}(u)$ will increase exponentially with $|u|$ [because of the exponential decay of $\phi(u)$], and we have to regularize the problem by cutting high frequencies. An appropriate weighting scheme would allow to take both effects into account.

6.2. *Truncation.* First, we truncate $\widetilde{\phi}$ to avoid the logarithm's explosion. Let

$$\widetilde{\mathcal{Y}}(u) := \log(-\log(T_{\omega_-,\omega_+}[|\widetilde{\phi}|^2](u))), \qquad u \in \mathbb{R} \setminus \{0\},$$

where the truncation operator $T_{\omega_-,\omega_+}$ with truncation levels $0 < \omega_- \leq \omega_+ < 1$ is defined via

$$T_{\omega_-,\omega_+}[f](u) = \begin{cases} \omega_+, & f(u) > \omega_+, \\ f(u), & \omega_- \leq f(u) \leq \omega_+, \\ \omega_-, & f(u) < \omega_-, \end{cases}$$

for any real-valued function $f$.

6.3. *Linearization.* Truncation allows us to linearize the problem. Set

$$\omega_\pm^*(u) := |\phi(u)|^2 \left(1 \pm \frac{2|\log|\phi(u)||}{1 + 2|\log|\phi(u)||}\right).$$

The following lemma holds:

LEMMA 6.1. *For any $u \in \mathbb{R} \setminus \{0\}$ and any $\omega_-(u), \omega_+(u)$ satisfying*

$$0 < \omega_- \leq \omega_-^* \leq \omega_+^* \leq \omega_+ < 1,$$

*the following inequality holds with probability one:*

$$|\widetilde{\mathcal{Y}}(u) - \mathcal{Y}(u) - \zeta_1(u)(|\widetilde{\phi}(u)|^2 - |\phi(u)|^2)| \leq \zeta_2(u)(|\widetilde{\phi}(u)|^2 - |\phi(u)|^2)^2,$$

*where*

$$\zeta_1(u) = 2^{-1}|\phi(u)|^{-2} \log^{-1}(|\phi(u)|)$$

*and*

$$\zeta_2(u) = 2 \max_{\xi \in \{\omega_-(u), \omega_+(u)\}} \left[\frac{1 + |\log(\xi)|}{\xi^2 \log^2(\xi)}\right].$$

Using the notation

$$\Delta(u) := |\widetilde{\phi}(u)|^2 - |\phi(u)|^2,$$

Lemma 6.1 can be reformulated as follows:

COROLLARY 6.2. *For any $u \in \mathbb{R} \setminus \{0\}$,*

$$\widetilde{\mathcal{Y}}(u) - \mathcal{Y}(u) = \zeta_1(u)\Delta(u) + Q(u), \tag{6.3}$$

*where*

$$|Q(u)| \leq \zeta_2(u)\Delta^2(u) \tag{6.4}$$

*with probability one.*

REMARK 6.1. Since $\phi(0)=1$ and $\phi(u)\to 0$ as $|u|\to\infty$, the behavior of truncation levels $\omega_-(u)$ and $\omega_+(u)$ in the vicinity of points $u=0$ and $u=\infty$ becomes important for determining the behavior of $\widetilde{\mathcal{Y}}(u)$ around these points. However, the values of $\widetilde{\mathcal{Y}}(u)$ around 0 will be discarded while estimating $\alpha$, and hence we do not need to know $\omega_+(u)$ for small $|u|$. As to $\omega_-(u)$ and $\omega_+(u)$ for large $u$, they can be constructed if some prior information on the Blumenthal–Getoor index $\alpha$ and the function $\pi(u)=\eta\tau(u)$ is available. For instance, if $0<\underline{\alpha}\le\alpha\le\overline{\alpha}\le 2$ and $0<\pi_-\le \mathrm{Re}[\pi(u)]\le\pi_+$ for all $|u|>u_0$ with large enough $u_0>0$, then one can take

$$\omega_-(u)=C_1 e^{-2\pi_+|u|^{\overline{\alpha}}}|u|^{-\overline{\alpha}},\qquad |u|>u_0,$$
$$\omega_+(u)=C_2 e^{-2\pi_-|u|^{\underline{\alpha}}},\qquad |u|>u_0,$$

with some constants $C_1>0$ and $C_2$ depending on $\pi_+$ and $\pi_-$, respectively. While a prior upper estimate $\overline{\alpha}$ for $\alpha$ appears also in the minimax rates of convergence proved in Section 6.6, a lower estimate $\underline{\alpha}$ turns out to be irrelevant for the convergence rates.

Note that the slope coefficient $\zeta_1$ grows exponentially with $|u|$. This means that the variance of $\widetilde{\mathcal{Y}}(u)$ grows exponentially as well and the values of $\widetilde{\mathcal{Y}}(u)$ with large $|u|$ and should be discarded when estimating $\alpha$.

6.4. *Spectral cut-off estimation.* Taking into account the special semi-linear structure of (6.2) together with a heteroscedastic variance of $\widetilde{\mathcal{Y}}(u)$, we apply a weighted least squares method to estimate $\alpha$. Let $w^1(u)$ be a function supported on $[\epsilon,1]$ with some $\epsilon>0$ that satisfies

$$\int_0^1 w^1(u)\log(u)\,du=1,\qquad \int_0^1 w^1(u)\,du=0. \tag{6.5}$$

For any $U>0$ put

$$w^U(u)=U^{-1}w^1(uU^{-1})$$

and define an estimate $\widetilde{\alpha}_U$ of $\alpha$ as

$$\widetilde{\alpha}_U=\int_0^\infty w^U(u)\widetilde{\mathcal{Y}}(u)\,du. \tag{6.6}$$

It is instructive to see what happens with $\widetilde{\alpha}_U$ in the case of exact data, that is, $\widetilde{\mathcal{Y}}=\mathcal{Y}$. One can see that in this case the following decomposition holds:

$$\widetilde{\alpha}_U=\log(2\eta)\underbrace{\int_0^\infty w^U(u)\,du}_{0}+\alpha\underbrace{\int_0^\infty w^U(u)\log(u)\,du}_{1}+R_U$$

with

$$R_U := \int_0^\infty w^U(u) \log(\operatorname{Re} \tau(u))\, du. \tag{6.7}$$

So, even in the case of perfect observations we still have the “bias” term $R_U$ induced by model misspecification. Indeed, when applying the least squares method we ignore a nonlinearity caused by $R_U$ and treat the problem as being linear. This is, however, only justified if $R_U$ is small. In fact, $R_U$ can be made small by taking large values of $U$.

6.5. *Further specification of the model class.* In order to rigorously study the complexity of the underlying estimation problem, we have to make further assumptions about the model class. Let us consider a class of Lévy models $\mathcal{A}(\overline{\alpha}, \eta_-, \eta_+, \varkappa)$ with

$$\psi(u) = \mathrm{i}\mu u + \vartheta(u), \qquad \vartheta(u) = -\eta|u|^\alpha \tau(u), \qquad u \in \mathbb{R}, \tag{6.8}$$

where $0 < \alpha \le \overline{\alpha} \le 2$,

$$0 < \eta_- \le \eta \le \eta_+ < \infty \tag{6.9}$$

and

$$|1 - \tau(u)| \lesssim \frac{1}{|u|^\varkappa}, \qquad |u| \to \infty, \tag{6.10}$$

for some $0 < \varkappa \le \alpha$. We will write

$$(\alpha, \eta, \tau) \in \mathcal{A}(\overline{\alpha}, \eta_-, \eta_+, \varkappa)$$

to indicate that the Lévy process with the characteristics $(\alpha, \eta, \tau)$ is in the class $\mathcal{A}$. The following proposition shows that conditions (6.8), (6.9) and (6.10) can be in fact rephrased in terms of the Lévy density of a $\mathcal{A}(\overline{\alpha}, \eta_-, \eta_+, \varkappa)$ process.

Proposition 6.3. *Let $\nu(x)$ be the Lévy density of a Lévy process satisfying (6.8) where the function $\tau$ fulfills*

$$\tau(u) = 1 + D_\pm u^{-\kappa} + o(|u|^{-\kappa}), \qquad u \to \pm\infty, \tag{6.11}$$

*with some constants $D_+$ and $D_-$. Then*

$$\int_{|x|<\epsilon} x^2 \nu(x)\, dx = c\epsilon^{2-\alpha} \theta(\epsilon), \tag{6.12}$$

*where $c > 0$ is a constant depending on $\eta$ and $\alpha$ and the function $\theta(\epsilon)$ satisfies*

$$|\theta(\epsilon) - 1| \lesssim |\epsilon|^\kappa, \qquad \epsilon \to 0.$$

As will be shown in the next two sections, even in the class $\mathcal{A}(\overline{\alpha}, \eta_-, \eta_+, \varkappa)$ the problem of estimating $\alpha$ is severely ill-posed, that is, a small perturbation $\varepsilon$ in data may lead (in worst case) to $\log^{-\varkappa/\overline{\alpha}}(1/\varepsilon)$ distance between $\alpha$ and its best estimate. On other hand, it turns out that our estimate $\widetilde{\alpha}_U$ achieves the best possible rates of convergence in the class $\mathcal{A}(\overline{\alpha}, \eta_-, \eta_+, \varkappa)$.

6.6. *Upper bounds.* Let us define

$$\varepsilon := \begin{cases} n^{-1}, & \text{under } \mathbb{P}, \\ \|\delta\|^2 + \sum_{j=1}^{n} \delta_j^2 \widetilde{\sigma}^2(y_j), & \text{under } \mathbb{Q}, \end{cases}$$

where $\|\delta\|^2 = \sum_{j=1}^n \delta_j^2$, $\widetilde{\sigma}(y_j) = e^{-y_j}\sigma(y_j)$ and $\delta_j = y_j - y_{j-1}$. In the case of calibration $\varepsilon$ comprises the level of the numerical interpolation error and of the statistical error simultaneously. In this section we will study the asymptotic behavior of the estimate $\widetilde{\alpha}_U = \widetilde{\alpha}_U(\varepsilon)$ defined in (6.6) as $\varepsilon \to 0$, $A := \min\{-y_0, y_n\} \to \infty$ and $e^{-A} \lesssim \|\delta\|^2$. Thus, it is assumed that the number of historical observations as well as the number of available options tend to infinity. First, we present an upper bound showing that our estimate $\widetilde{\alpha}_U$ with the "optimal" choice of the cut-off parameter $U$ converges to $\alpha$ with a logarithmic rate in $\varepsilon$.

THEOREM 6.4. *For $U = \overline{U}$ with*

$$\overline{U} = \left[\frac{1}{2\eta_+} \log(\varepsilon^{-1} \log^{-\beta}(1/\varepsilon))\right]^{1/\overline{\alpha}}$$

*and*

$$\beta = \begin{cases} 1 + \varkappa/\overline{\alpha}, & \text{under } \mathbb{P}, \\ (\varkappa + 4)/\overline{\alpha} - 1, & \text{under } \mathbb{Q}, \end{cases}$$

*it holds*

$$\sup_{(\alpha,\eta,\tau)\in\mathcal{A}(\overline{\alpha},\eta_-,\eta_+,\varkappa)} \mathrm{E}|\widetilde{\alpha}_{\overline{U}} - \alpha|^2 \lesssim \mathcal{R}(\varepsilon), \qquad \varepsilon \to 0, \tag{6.13}$$

*where*

$$\mathcal{R}(\varepsilon) = \left[\frac{1}{2\eta_+} \log \varepsilon^{-1}\right]^{-2\varkappa/\overline{\alpha}}.$$

REMARK 6.2. Since the rates are logarithmic it is usual to call the underlying estimation problem severely ill-posed. From a practical point of view, severely ill-posedness means that more observations are needed to reach the desired level of accuracy than for well-posed problems.

REMARK 6.3. As can be easily seen the convergence rates depend on $\overline{\alpha}$, a prior upper bound for $\alpha$. If there is no prior information on $\overline{\alpha}$ one may take $\overline{\alpha} = 2$.

REMARK 6.4. For symmetric stable processes we have $\tau(u) \equiv 1$ and it can be shown that the rates are parametric in this case, that is,

$$\sup_{(\alpha,\eta,\tau)\in\mathcal{A}(\overline{\alpha},\eta_-,\eta_+,\infty)} \mathrm{E}|\widetilde{\alpha}_{\overline{U}} - \alpha|^2 \lesssim \varepsilon, \qquad \varepsilon \to 0,$$

for some $\overline{U}$ depending on $\varepsilon$.

6.7. *Lower bounds.* Now we show that the rates obtained in the previous section are the best ones in the minimax sense for the class $\mathcal{A}(\overline{\alpha}, \eta_-, \eta_+, \varkappa)$.

THEOREM 6.5. *It holds*

$$\lim_{s\to 0} \liminf_{\varepsilon\to 0} \inf_{\widetilde{\alpha}} \sup_{(\alpha,\eta,\tau)\in\mathcal{A}(\overline{\alpha},\eta_-,\eta_+,\varkappa)} \delta_{n,s}^{-2}(\varepsilon)\mathrm{E}(|\widetilde{\alpha} - \alpha|^2) = O(1), \tag{6.14}$$

*where*

$$\delta_{n,s}(\varepsilon) = \left[\frac{1}{2\eta_+} \log \varepsilon^{-1}\right]^{-\varkappa/(\overline{\alpha}-s)}$$

*and the infimum is taken over all estimators $\widetilde{\alpha}$ of $\alpha$.*

6.8. *Asymptotic behavior.* In this section we complete the investigation of asymptotic properties of the estimate $\widetilde{\alpha}$ by proving its asymptotic normality. In the case of estimation under $\mathbb{P}$ we have the following:

THEOREM 6.6. *Denote*

$$\varsigma(\varepsilon, U) = \left[\varepsilon \int_0^\infty w^U(u) w^U(v) \zeta_1(u) \zeta_1(v) S(u,v)\, du\, dv\right]^{1/2}$$

*with*

$$\begin{aligned} S(u,v) := {} & \operatorname{Re}\phi(u-v) + \operatorname{Im}\phi(u+v) \\ & - (\operatorname{Re}\phi(u) + \operatorname{Im}\phi(u))(\operatorname{Re}\phi(v) + \operatorname{Im}\phi(v)). \end{aligned}$$

*Let $U = U(\varepsilon)$ be a sequence of cutoffs such that $\varsigma^{-1}(\varepsilon, U(\varepsilon))R_{U(\varepsilon)} \to 0$ as $\varepsilon \to 0$. Then*

$$\varsigma^{-1}(\varepsilon, U(\varepsilon))(\widetilde{\alpha}_{U(\varepsilon)} - \alpha) \sim \mathcal{N}(0,1), \qquad \varepsilon \to 0.$$

REMARK 6.5. The choice of $U(\varepsilon)$ is based on the following reasoning. On the one hand, we have to require that $\varsigma^{-1}(\varepsilon, U)R_U \to 0$ in order to ensure that $\varsigma^{-1}(\varepsilon, U)(\widetilde{\alpha}_U - \alpha)$ has asymptotically zero expectation. On the other hand, the variance of $\varsigma^{-1}(\varepsilon, U)\widetilde{\alpha}_U$ should converge as $\varepsilon \to 0$, and the limit must be bounded and nondegenerated.

Remark 6.6. Given an estimate $\widetilde{\phi}$ of $\phi$ and some $U = U(\varepsilon)$ such that $|\widetilde{\phi}(u)| \neq 0$ on $[-U, U]$ and $|\widetilde{\phi}(u)| \neq 1$ on $[-U, U] \setminus \{0\}$, we can estimate the norming factor $\varsigma(\varepsilon, U)$ for $\widetilde{\alpha}_U$ via

$$\varsigma(\varepsilon, U) := \left[\varepsilon \int_0^\infty w^U(u) w^U(v) \widetilde{\zeta}_1(u) \widetilde{\zeta}_1(v) \widetilde{S}(u, v)\, du\, dv\right]^{1/2}$$

with

$$\begin{aligned}\widetilde{S}(u, v) := {} & \operatorname{Re} \widetilde{\phi}(u - v) + \operatorname{Im} \widetilde{\phi}(u + v) \\ & - (\operatorname{Re} \widetilde{\phi}(u) + \operatorname{Im} \widetilde{\phi}(u))(\operatorname{Re} \widetilde{\phi}(v) + \operatorname{Im} \widetilde{\phi}(v))\end{aligned}$$

and

$$\widetilde{\zeta}_1(u) := |\widetilde{\phi}(u)|^{-2} \log^{-1}(|\widetilde{\phi}(u)|^2).$$

A similar result can be proved in the case of calibration as well.

6.9. *Processes with a nonzero diffusion part.* In fact, spectral calibration algorithm allows us to treat more general models with a nonzero diffusion part. Let $\mathcal{A}(\overline{a}, \overline{\alpha}, \eta_-, \eta_+, \varkappa)$ be a class of Lévy processes with the characteristic exponent of the form

$$\psi_a(u) = \mathrm{i}\mu u - a^2 u^2/2 + \vartheta(u), \qquad \vartheta(u) = -\eta |u|^\alpha \tau(u), \qquad u \in \mathbb{R}, \tag{6.15}$$

where $0 < a < \overline{a}$ and conditions (6.9) and (6.10) are fulfilled. We will write $(a, \alpha, \eta, \tau) \in \mathcal{A}(\overline{a}, \overline{\alpha}, \eta_-, \eta_+, \varkappa)$ to indicate that a Lévy process with the characteristic exponent (6.15) belongs to $\mathcal{A}(\overline{a}, \overline{\alpha}, \eta_-, \eta_+, \varkappa)$.

Assume first that $\phi_a(u) = \exp(\psi_a(u))$ is known exactly. Define

$$\mathcal{L}(u) := \log(|\phi_a(u)|^2) = -a^2 u^2 + 2 \operatorname{Re}[\vartheta(u)]$$

and

$$\mathcal{L}_\xi(u) := \xi^2 \mathcal{L}(u) - \mathcal{L}(\xi u) = \log(|\phi_a(u)|^{2\xi^2} / |\phi_a(\xi u)|^2) =: \log(\rho_\xi(u))$$

for some $\xi > 1$. It obviously holds

$$\mathcal{L}_\xi(u) = -\eta |u|^\alpha (\xi^2 \operatorname{Re}[\tau(u)] - \xi^\alpha \operatorname{Re}[\tau(\xi u)]) = -c_\xi(\alpha) |u|^\alpha \tau_\xi(u),$$

where $c_\xi(\alpha) = \eta(\xi^2 - \xi^\alpha)$, and $\tau_\xi(u)$ fulfills

$$|1 - \tau_\xi(u)| \lesssim \frac{1}{|u|^\varkappa}, \qquad |u| \to \infty. \tag{6.16}$$

Thus, $\mathcal{L}_\xi(u)$ has a structure similar to the structure of $\vartheta(u)$ in (6.8) and we can carry over the results of the previous section to a more general model (6.15) by defining

$$\widetilde{\mathcal{Y}}_\xi(u) := \log(-\log(T_{\omega_-, \omega_+}[\widetilde{\rho}_\xi](u))),$$

where $\widetilde{\rho}_\xi(u) = |\widetilde{\phi}(u)|^{2\xi^2}/|\widetilde{\phi}(\xi u)|^2$ with $\widetilde{\phi}$ being an estimate of $\phi_a$. Define

$$\widetilde{\alpha}_{\xi,U} = \int_0^\infty w^U(u)\widetilde{\mathcal{Y}}_\xi(u)\,du. \tag{6.17}$$

The following two theorems are extensions of Theorems 6.4 and 6.5, respectively, to the case of Lévy models with a nonzero diffusion part.

THEOREM 6.7. *For $U = \overline{U}$ with*

$$\overline{U} = \left[\frac{1}{2\overline{a}}\log(\varepsilon^{-1}\log^{-\beta}(1/\varepsilon))\right]^{1/2}$$

*and $\beta = 1 + \varkappa/2$, it holds*

$$\sup_{(a,\alpha,\eta,\tau)\in\mathcal{A}(\overline{a},\overline{\alpha},\eta_-,\eta_+,\varkappa)} \mathrm{E}|\widetilde{\alpha}_{\xi,\overline{U}} - \alpha|^2 \lesssim \mathcal{R}(\varepsilon), \qquad \varepsilon \to 0, \tag{6.18}$$

*where*

$$\mathcal{R}(\varepsilon) = c_\xi^{-1}(\overline{\alpha})\left[\frac{1}{2\overline{a}}\log\varepsilon^{-1}\right]^{-\varkappa}.$$

THEOREM 6.8. *It holds*

$$\liminf_{\varepsilon\to 0}\inf_{\widetilde{\alpha}} \sup_{(a,\alpha,\eta,\tau)\in\mathcal{A}(\overline{a},\overline{\alpha},\eta_-,\eta_+,\varkappa)} \delta_n^{-2}(\varepsilon)\mathrm{E}(|\widetilde{\alpha} - \alpha|^2) = O(1), \tag{6.19}$$

*where*

$$\delta_n(\varepsilon) = \left[\frac{1}{2\overline{a}}\log\varepsilon^{-1}\right]^{-\varkappa/2},$$

*and the infimum is taken over all estimators $\widetilde{\alpha}$ of $\alpha$.*

As can be seen, the estimate $\widetilde{\alpha}_{\xi,\overline{U}}$ is consistent as long as $\overline{\alpha} < 2$. The nearer is $\overline{\alpha}$ to 2, the closer is the constant $c_\xi(\overline{\alpha})$ to zero and the more difficult becomes the estimation problem.

**7. Adaptive procedure.** Minimax results obtained in the previous sections show the complexity of the underlying estimation problem but are not very helpful in practice. Putting aside the fact that they are related to the performance of the procedure in the worst situation (worst case scenario) which is not necessarily the case for the given model from $\mathcal{A}(\overline{\alpha}, \eta_-, \eta_+, \varkappa)$, the choice of $U$ suggested there depends on $\overline{\alpha}$, is asymptotic and likely to be inefficient for small sample sizes. In this section we propose an adaptive procedure for choosing the cut-off parameter $U$. First, let us fix a sequence of cut-off parameters $U_1 > U_2 > \cdots > U_K$ and define

$$\widetilde{\alpha}_k = \int_0^\infty w^{U_k}(u)\widetilde{\mathcal{Y}}(u)\,du, \qquad k = 1, \ldots, K.$$

We suggest a method based on the combination of multiple testing and aggregation ideas [see Belomestny and Spokoiny (2007)]. Namely, for the sequence of estimates $\widetilde{\alpha}_k$ consider a sequence of nested hypothesis $H_k : \alpha_1 = \cdots = \alpha_k = \alpha$ where

$$\alpha_k = \int_0^\infty w^{U_k}(u)\mathcal{Y}(u)\,du, \qquad k = 1, \ldots, K.$$

The hypothesis $H_k$ basically means that $R_{U_i} = 0$ for $i = 1, \ldots, k$. The procedure is sequential; we put $\widehat{\alpha}_1 = \widetilde{\alpha}_1$ and start with $k = 2$ and at each step $k$ the hypothesis $H_k$ is tested against $H_{k-1}$. For testing $H_k$ against $H_{k-1}$ we check that the previously constructed adaptive estimate; $\widehat{\alpha}_{k-1}$ belongs to the confidence intervals built on $\widetilde{\alpha}_k$. Then we put

$$\widehat{\alpha}_k = \gamma_k \widetilde{\alpha}_k + (1 - \gamma_k)\widehat{\alpha}_{k-1}. \tag{7.1}$$

The mixing parameter $\gamma_k$ is defined using a measure of statistical difference between $\widehat{\alpha}_{k-1}$ and $\widetilde{\alpha}_k$

$$\gamma_k := \mathcal{K}(T_k/\mathcal{V}_k), \qquad T_k := (\widetilde{\alpha}_k - \widehat{\alpha}_{k-1})^2/\sigma_k^2,$$

where $\sigma_k^2$ is the variance of $\widetilde{\alpha}_k$, $\mathcal{K}$ is a kernel supported on $[0,1]$ and $\{\mathcal{V}_k\}$ is a set of critical values. In particular, $\gamma_k$ is equal to zero if $H_k$ is rejected; that is, $\widehat{\alpha}_{k-1}$ lies outside the confidence interval around $\widetilde{\alpha}_k$. The final estimate is equal to $\widehat{\alpha}_K$.

7.1. *Choice of the critical values* $\mathcal{V}_k$. The critical values $\mathcal{V}_1, \ldots, \mathcal{V}_{K-1}$ are selected by a reasoning similar to the standard approach of hypothesis testing theory. We would like to provide prescribed performance of the procedure under the simplest (null) hypothesis. In the considered set-up, the null means that

$$\alpha_1 = \cdots = \alpha_K = \alpha. \tag{7.2}$$

In this case it is natural to expect that the estimate $\widehat{\alpha}_k$ coming out of the first steps of the procedure until index $k$ is close to the nonadaptive counterpart $\widetilde{\alpha}_k$.

To give a precise definition we need to specify a loss function. Suppose that the risk of estimation for an estimate $\widehat{\alpha}$ of $\alpha$ is measured by $\mathrm{E}|\widehat{\alpha} - \alpha|^{2r}$ for some $r > 0$. It is not difficult to show that under the null hypothesis (7.2), each estimate $\widetilde{\alpha}_k$ asymptotically fulfills

$$\varepsilon^{-1/2}(\widetilde{\alpha}_k - \alpha) \sim \mathcal{N}(0, \sigma_k^2), \qquad \varepsilon \to 0.$$

For example, in the case of estimation under $\mathbb{P}$ one can prove (see the proof of Proposition 6.6) that

$$\sigma_k^2 = \int_0^\infty \int_0^\infty w^{U_k}(u)\zeta_1(u) w^{U_k}(v)\zeta_1(v) S(u,v)\,du\,dv \tag{7.3}$$

with

$$S(u,v) := \operatorname{Re}\phi(u-v) + \operatorname{Im}\phi(u+v) \\ - (\operatorname{Re}\phi(u) + \operatorname{Im}\phi(u))(\operatorname{Re}\phi(v) + \operatorname{Im}\phi(v)).$$

Therefore,

$$\mathrm{E}_0|\sigma_{k,\varepsilon}^{-2}(\widetilde{\alpha}_k - \alpha)^2|^r \approx C_r,$$

where $\sigma_{k,\varepsilon}^2 = \varepsilon\sigma_k^2$, $C_r = E|\xi|^{2r}$, and $\xi$ is the standard normal. We require the parameters $\mathcal{V}_1, \dots, \mathcal{V}_{K-1}$ of the procedure to satisfy

$$\mathrm{E}_0|\sigma_{k,\varepsilon}^{-2}(\widehat{\alpha}_k - \widetilde{\alpha}_k)^2|^r \le \gamma C_r, \qquad k = 2, \dots, K. \tag{7.4}$$

Here $\gamma$ stands for a preselected constant having the meaning of a confidence level of the procedure. This gives us $K-1$ conditions to fix $K-1$ parameters.

Our definition still involves two parameters $\gamma$ and $r$. It is important to mention that their choice is subjective and there is no way for an automatic selection. A proper choice of the power $r$ for the loss function as well as the “confidence level” $\gamma$ depends on the particular application and on the additional subjective requirements of the procedure.

## 8. Simulations.

8.1. *Estimation of the fractional order from a time series.* Let us consider the generalized hyperbolic (GH) Lévy model which was introduced in a series of papers [Eberlein and Keller (1995), Eberlein, Keller and Prause (1998) and Eberlein and Prause (2002)] and emerged from extensive empirical investigations of financial time series. See also Eberlein (2001) for a survey on a number of analytical aspects of this model. The characteristic function $\Phi_{\mathrm{GH}}$ of increments in the GH Lévy model with parameters $(\kappa, \beta, \delta, \lambda)$ is given by

$$\Phi_{\mathrm{GH}}(u) = e^{\mathrm{i}\mu u} \frac{(\sqrt{\kappa^2 + \beta^2})^\lambda}{(\sqrt{\kappa^2 - (\beta + \mathrm{i}u)^2})^\lambda} \frac{K_\lambda(\delta\sqrt{\kappa^2 - (\beta + \mathrm{i}u)^2})}{K_\lambda(\delta\sqrt{\kappa^2 + \beta^2})},$$

where $K$ is the modified bessel function of the second kind. $\Phi_{\mathrm{GH}}$ has the Lévy–Khintchine representation of the form,

$$\Phi_{\mathrm{GH}}(u) = \exp\left(ibu + \int_{-\infty}^{\infty} (e^{\mathrm{i}ux} - 1 - \mathrm{i}ux) g(x)\, dx\right).$$

Note that this model does not contain a Gaussian component $a^2u^2/2$. Function $g(x)$, the density of the corresponding Lévy measure, can be represented [see Eberlein (2001)] in an integral form. From this representation the following expansion for $\rho(x) = x^2 g(x)$ can be obtained;

$$\rho(x) = \frac{\delta}{\pi} + \frac{\lambda + 1/2}{2}|x| + \frac{\delta\beta}{\pi}x + o(|x|), \qquad x \to 0.$$

A direct consequence of this expansion is that

$$\int_{|x|>\varepsilon} g(x)\,dx \asymp 1/\varepsilon, \qquad \varepsilon \to 0,$$

and hence the fractional order of the GH Lévy model is equal to 1. In our simulation study we simulate GH Lévy process $X$ with $\beta = 0$, $\lambda = 1$ and different pairs of $\kappa$ and $\delta$ at $n+1$ equidistant points $\{0, \Delta, \ldots, n\Delta\}$. Upon that we construct the empirical characteristic function of increments;

$$\widetilde{\phi}(u) = \frac{1}{n}\sum_{k=1}^{n} e^{\mathrm{i}u(X_{k\Delta} - X_{(k-1)\Delta})}.$$

Following the description of the spectral estimation algorithm, define

$$\widetilde{\mathcal{Y}}(u) := \log(-\log(T_{\omega_-,\omega_+}[|\widetilde{\phi}|^2](u))),$$

where truncation levels $\omega_-$, and $\omega_+$ are constant in $u$ and are equal to 0.01 and 0.95, respectively. In fact, for practical applications with a medium sample sizes $n$, the choice of these levels is not crucial. Now consider the following minimization problem:

$$(l_0^U, l_1^U) = \arg\min_{l_0, l_1} \int_0^U \overline{w}^U(u)(\widetilde{\mathcal{Y}}(u) - l_1 \log(u) - l_0)^2\,du, \tag{8.1}$$

where $\overline{w}^U(u) = U^{-1}\overline{w}^1(U^{-1}u)$, and $\overline{w}^1(u) = u\mathbf{1}_{\{\epsilon \le u \le 1\}}$ for some $\epsilon > 0$. An estimate for the fractional order is defined as $\widetilde{\alpha}^U = l_1^U$. It is not difficult to show that $\widetilde{\alpha}^U$ is of the form,

$$\widetilde{\alpha}^U = \int_0^\infty w^U(u)\widetilde{\mathcal{Y}}(u)\,du$$

with $w^U(u) = U^{-1}w^1(U^{-1}u)$ and $w^1(u) = \overline{w}^1(u)[A_1\log(u) - A_2]$ where $A_1$ and $A_2$ are two positive constants such that $w^1(u)$ satisfies conditions (6.5). Let $U_1 > U_2 > \cdots > U_K$ be an exponentially decreasing sequence of cut-offs and $\widetilde{\alpha}_1, \ldots, \widetilde{\alpha}_K$ be the corresponding sequence of estimates. Following (7.1), we construct a sequence of aggregated estimates $\widehat{\alpha}_1, \ldots, \widehat{\alpha}_K$ using a triangle kernel and a set of critical values $\mathcal{V}_1, \ldots, \mathcal{V}_K$ computed by (7.4). The variances $\{\sigma_k^2\}$ in (7.3) are estimated from above using a bound for $\zeta_1$. Box plots of $\widehat{\alpha} = \widehat{\alpha}_K$ based on 500 trials for different $n$ and different pairs of $\kappa$ and $\delta$ are shown in Figure 1.

8.2. *Estimation of the fractional order from options data.* In the case of calibration (estimation under $\mathbb{Q}$) we compute first the prices of $n$ call options,

$$\mathcal{C}(y_k, T) = S\mathrm{E}^{\mathbb{Q}}[(e^{Y_T} - e^{y_k})^+], \qquad k = 1, \ldots, n,$$

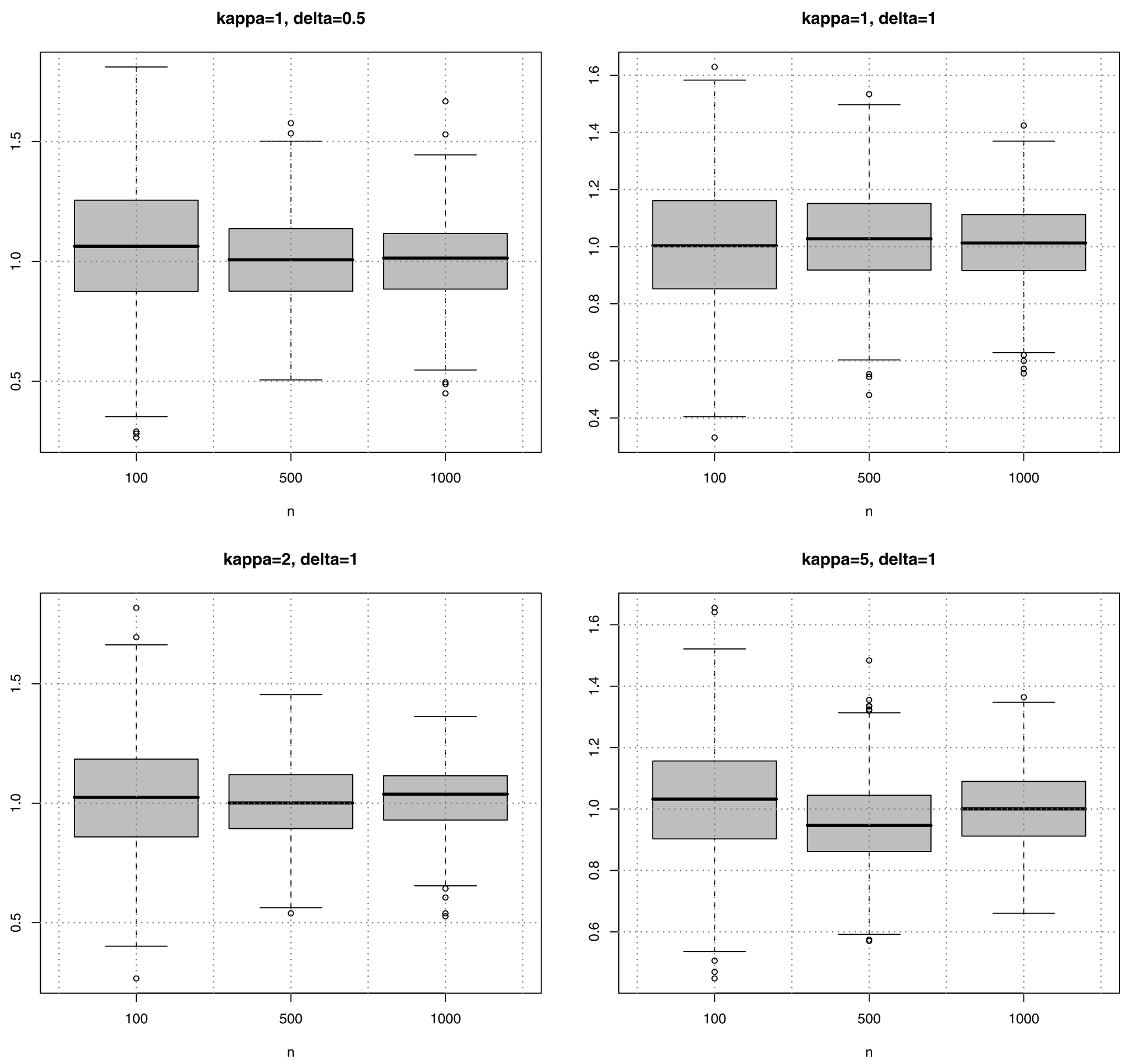

FIG. 1. *Box plots of the estimate $\widehat{\alpha}$ under $\mathbb{P}$ for different sample sizes and different parameters of the GH process.*

using formula (3.3) where the underlying process $Y$ follows a GH Lévy model (parameters will be specified later on): $S = 1$, $T = 0.25$ and $r = 0.06$. The log-moneyness design $(y_i)$ is chosen to be normally distributed with zero mean and variance $1/3$ and reflects the structure of the option market where much more contracts are settled at the money than in or out of money. Finally, we simulate

$$\mathcal{O}_T(y_j) = O_T(y_j) + \sigma(y_j)\xi_j, \qquad j = 1, \dots, n,$$

where $\xi_j$ are standard normal, $O_T$ is defined in (3.2) and $\sigma(y) = [\overline{\sigma} O_T(y)]^2$.

In the first step of our estimation procedure we find the function $\widehat{O}$ among all functions $O$ with two continuous derivatives as the minimizer of the penalized residual sum of squares

$$\mathrm{RSS}(O, L) = \sum_{i=0}^{n+1} (\mathcal{O}_T(y_i) - O(y_i))^2 + L \int_{y_0}^{y_{n+1}} [O''(u)]^2 \, du, \tag{8.2}$$

where $y_0 \ll y_1$ and $y_{n+1} \gg y_n$ are two extrapolated points with artificial values $O_{n+1} = O_0 = 0$. The first term in (8.2) measures closeness to the data, while the second term penalizes curvature in the function, and $L$ establishes a trade-off between the two. The two special cases are $L = 0$ when $\widehat{O}$ interpolates the data, and $L = \infty$ when a straight line using ordinary least squares is fitted. In our numerical example we use the R package `p-splines` with the choice of $L$ that minimizes the generalized cross-validation criterion. It can be shown that (8.2) has an explicit, finite dimensional, unique minimizer which is a natural cubic spline with knots at the values of $y_i, i = 1, \ldots, n$. Since the solution of (8.2) is a natural cubic spline, we can write

$$\widehat{O}(y) = \sum_{j=1}^{n} \theta_j \beta_j(y),$$

where $\beta_j(y), j = 1, \ldots, n$, is a set of basis functions representing the family of natural cubic splines. We estimate $\mathbf{F}[\widehat{O}](v + \mathrm{i})$ by

$$\mathbf{F}[\widehat{O}](v + \mathrm{i}) = \sum_{j=1}^{n} \theta_j \mathbf{F}[e^{-y}\beta_j(y)](v).$$

Although $\mathbf{F}[e^{-y}\beta_j(y)]$ can be computed in closed form, we just use the fast Fourier transform (FFT) and compute $\mathbf{F}[\widehat{O}](v + \mathrm{i})$ on a fine dyadic grid. On the same grid one can compute

$$\widetilde{\psi}(v) := \frac{1}{T} \log(1 + v(v + \mathrm{i})\mathbf{F}[\widehat{O}](v + \mathrm{i})), \qquad v \in \mathbb{R}, \tag{8.3}$$

where $\log(\cdot)$ is taken in such a way that $\widetilde{\psi}(v)$ is continuous with $\widetilde{\psi}(-\mathrm{i}) = 0$. Now we can follow the road map of the adaptive spectral calibration algorithm and get an estimate for the fractional order of the underlying GH Lévy model. In Figure 2 box plots of the final estimate $\widehat{\alpha} = \widehat{\alpha}_K$ based on 500 Monte Carlo trials are shown in the case of the underlying GH Lévy model with parameters $\beta = -1$, $\lambda = 1$ and different $\kappa, \delta$. Sample size $n$ is equal to 1000 and noise level $\overline{\sigma}$ takes values in the set $\{1, 10, 20\}$. The estimate $\widehat{\alpha}$ is obviously biased because of numerical errors (due to the approximation of Fourier integral and linearization).

8.3. *Processes with a nonzero diffusion part.* Turn now to the class of Lévy processes containing a nonzero diffusion part which was treated in Section 6.9. The only algorithmic difference to the case of processes with zero diffusion part is that now we first fix some $\xi > 1$ and compute

$$\widetilde{\mathcal{Y}}_\xi(u) := \log(-\log(T_{\omega_-,\omega_+}[|\widetilde{\rho}_\xi|^2](u)))$$

instead of $\widetilde{\mathcal{Y}}(u)$ where $\widetilde{\rho}_\xi(u) = |\widetilde{\phi}(u)|^{2\xi^2}/|\widetilde{\phi}(\xi u)|^2$ with $\widetilde{\phi}$ being an estimate of $\phi_a$. In the estimation procedure we consider only the set of $u$ with $|\widetilde{\phi}(\xi u)| > 0$. Note that this set is smaller than the set where $|\widetilde{\phi}(u)| > 0$ since $\xi > 1$. It is also intuitively clear that more observations are needed to estimate $\widetilde{\rho}_\xi$ with the same quality as $|\widetilde{\phi}(u)|^2$, and therefore the first problem is likely to be computationally more difficult. This conjecture is supported by

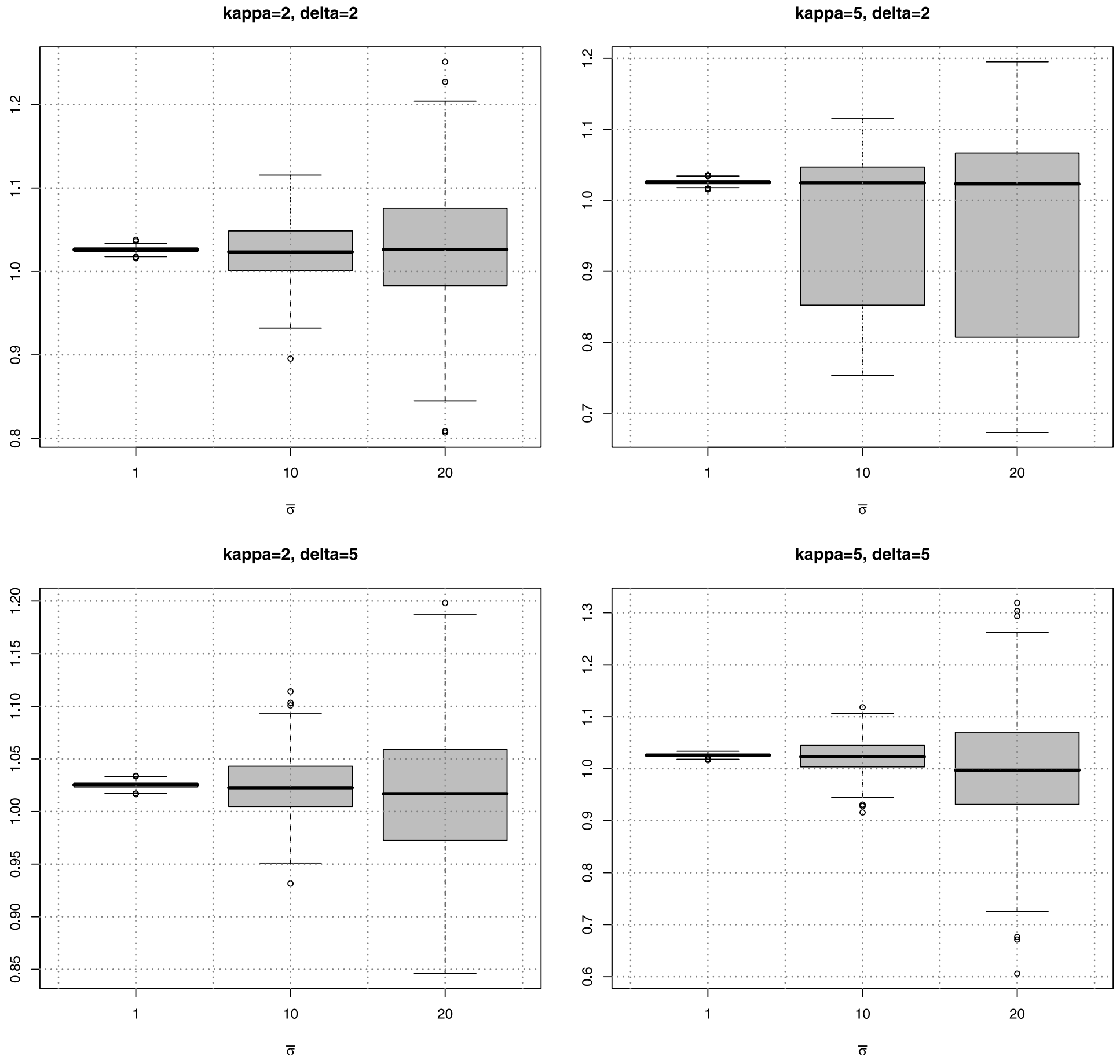

FIG. 2. *Box plots of the estimate $\widehat{\alpha}$ under $\mathbb{Q}$ for different noise levels and different sets of parameters of the underlying GH Lévy process.*

our simulation study as well. Figure 3 shows the boxplots of two estimates $\widehat{\alpha}$ and $\widehat{\alpha}_\xi$ based on 500 samples under historical measure $\mathbb{P}$ from the GH Lévy model with zero diffusion part (left) and with the diffusion parameter $a$ equal to 0.1 (right), remaining parameters $\lambda$, $\beta$, $\kappa$ and $\delta$ being equal to 1, 0, 1 and 4, respectively. The estimate $\widehat{\alpha}_\xi$ is constructed from the estimates $\widetilde{\alpha}_{U_1,\xi}, \ldots, \widetilde{\alpha}_{U_K,\xi}$ [we use $\xi = 2$ and $\overline{w}^1(u) = u\mathbf{1}_{\{\epsilon<u\leq 1\}}$ in (8.1)] via the stagewise aggregation procedure as described in Section 7. We took $K = 30$, $U_k = 100(1.25)^{-(k-1)}, k = 1, \ldots, K$ and $\mathcal{K}(x) = (1-x)1_{\{0\leq x\leq 1\}}$. As to the critical values, they are determined via (7.4) with $r = 1$, $\gamma = 0.5$. Note that while the difference between $\widehat{\alpha}_\xi$ and $\widehat{\alpha}$ is rather pronounced for small sample sizes, it almost disappears for sample sizes as large as 1000.

# APPENDIX

**A.1. Proof of Lemma 6.1.** For any positive $\omega_-$ and $\omega_+$ satisfying $\omega_-(u) \leq |\phi(u)|^2 \leq \omega_+(u)$, we have

$$
\begin{aligned}
&|\widetilde{\mathcal{Y}}(u) - \mathcal{Y}(u) - \zeta_1(u)(T_{\omega_-,\omega_+}[|\widetilde{\phi}|^2](u) - |\phi(u)|^2)| \\
&\quad \leq \frac{\zeta_2(u)}{2}(T_{\omega_-,\omega_+}[|\widetilde{\phi}|^2](u) - |\phi(u)|^2)^2 \\
&\quad \leq \frac{\zeta_2(u)}{2}(|\widetilde{\phi}(u)|^2 - |\phi(u)|^2)^2.
\end{aligned}
$$

Furthermore,

$$
||\widetilde{\phi}(u)|^2 - T_{\omega_-,\omega_+}[|\widetilde{\phi}|^2](u)| \leq ||\widetilde{\phi}(u)|^2 - |\phi(u)|^2|, \qquad u \in \mathbb{R}^d,
$$

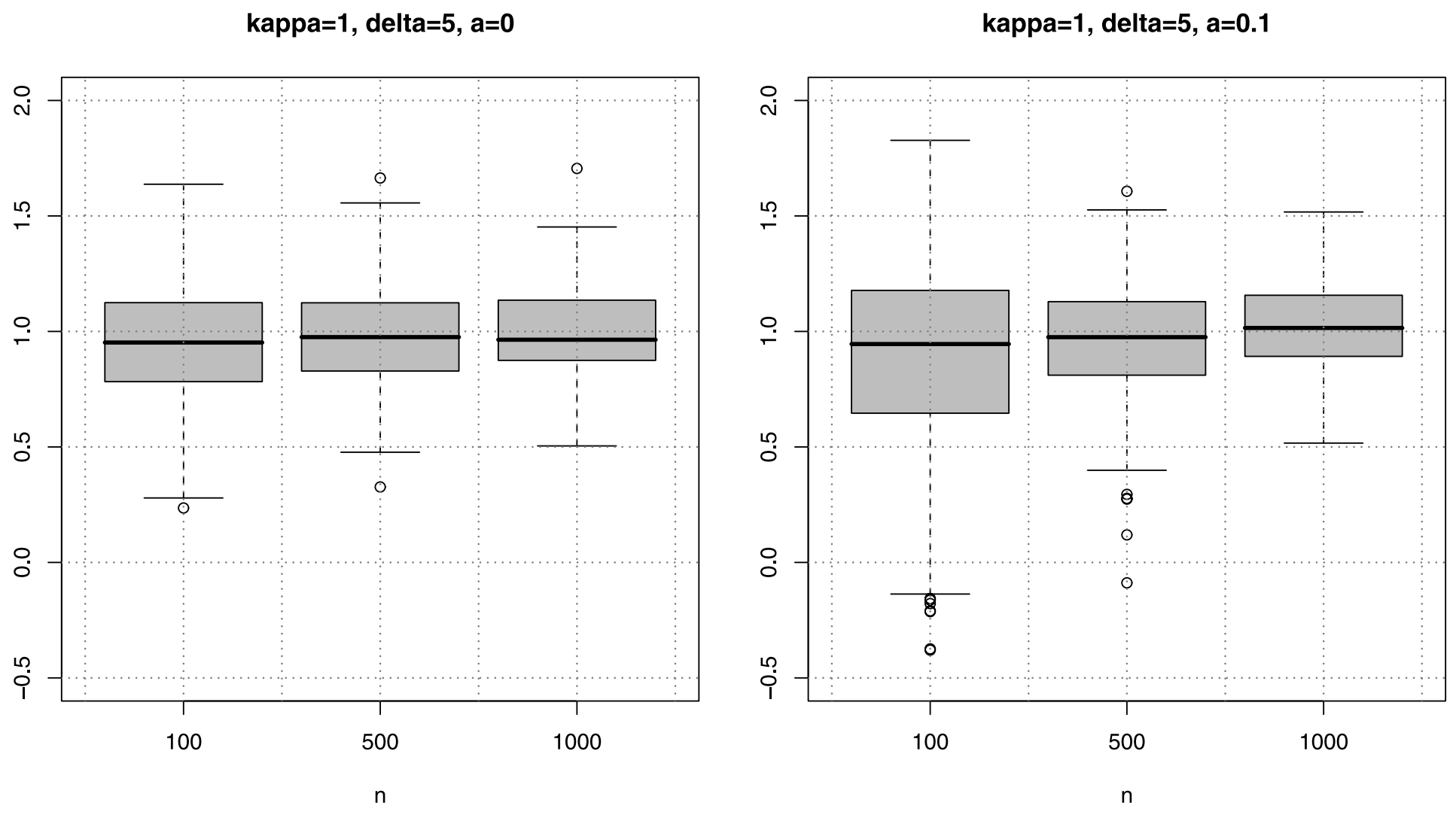

FIG. 3. *Box plots of the estimates $\widehat{\alpha}$ (left) and $\widehat{\alpha}_\xi$ (right) under $\mathbb{P}$ for different sample sizes $n$.*

and it holds on the set $|\widetilde{\phi}(u)|^2 \notin [\omega_-, \omega_+]$

$$||\widetilde{\phi}(u)|^2 - |\phi(u)|^2| \geq \min\{|\phi(u)|^2 - \omega_-, \omega_+ - |\phi(u)|^2\}.$$

Thus,

$$\zeta_1(u)||\widetilde{\phi}(u)|^2 - T_{\omega_-,\omega_+}[|\widetilde{\phi}|^2](u)| \leq \frac{\zeta_2(u)}{2}||\widetilde{\phi}(u)|^2 - |\phi(u)|^2|^2$$

on the set $|\widetilde{\phi}(u)|^2 \notin [\omega_-, \omega_+]$, provided that

$$2|\phi(u)|^2|\log(|\phi(u)|)|\min\{|\phi(u)|^2 - \omega_-, \omega_+ - |\phi(u)|^2\} \geq \frac{|\phi(u)|^4\log^2(|\phi(u)|^2)}{1 + |\log(|\phi(u)|^2)|},$$

that is,

$$\min\left\{1 - \frac{\omega_-}{|\phi(u)|^2}, \frac{\omega_+}{|\phi(u)|^2} - 1\right\} \geq \frac{\log(|\phi(u)|^2)}{1 + |\log(|\phi(u)|^2)|}.$$

**A.2. Proof of Proposition 6.3.** Without loss of generality we can assume that $\mu = 0$ in (6.8). Denote

$$\rho(x) = \left(1 - \frac{\sin x}{x}\right)\nu(x);$$

then $\rho$ is, up to a scaling factor, the density of some probability distribution with the characteristic function $\zeta_\rho\widetilde{\psi}(u)$ where $\zeta_\rho$ is a positive constant and

$$\widetilde{\psi}(u) = \int_{-1}^{1}(\psi(u) - \psi(u+w))\,dw.$$

Due to (6.11) the following asymptotic expansion holds:

$$\begin{aligned}\widetilde{\psi}(u) &= |u|^\alpha\tau(u)\int_{-1}^{1}\left[1 - \left|1 + \frac{w}{u}\right|^\alpha \frac{\tau(u+w)}{\tau(u)}\right]dw\\ &= C_\pm(\alpha,\kappa)|u|^{\alpha-2}[1 + O(|u|^{-\kappa})], \qquad u \to \pm\infty,\end{aligned}$$

with some constants $C_+$ and $C_-$ depending on $\alpha$ and $\kappa$. We consider separately two cases.

*Case* $0 < \alpha < 1$. Note that in this case $\widetilde{\psi}(u)$ is integrable on $\mathbb{R}$ and the Fourier inversion formula implies

$$\rho(x) = \frac{\zeta_\rho}{2\pi}\int_{-\infty}^{\infty}(\exp(-\mathrm{i}xu) - 1)\widetilde{\psi}(u)\,du$$

since $\rho(0)=0$. We have for any positive number $a$,

$$\begin{aligned}
&\int_{-\infty}^{\infty}(\exp(-\mathrm{i}xu)-1)\widetilde{\psi}(u)\,du\\
&\quad=\int_{|u|\le a}(\exp(-\mathrm{i}xu)-1)\widetilde{\psi}(u)\,du+\int_{|u|>a}(\exp(-\mathrm{i}xu)-1)\widetilde{\psi}(u)\,du\\
&\quad=: I_1+I_2,
\end{aligned}$$

where $|I_1|\lesssim|x|\lesssim|x|^{1-\alpha+\kappa}$ for $x\to 0$ provided that $\kappa\le\alpha$. Furthermore,

$$\begin{aligned}
I_2&=C_\pm(\alpha,\kappa)\int_{|u|>a}(\exp(-\mathrm{i}xu)-1)|u|^{\alpha-2}\,du+O(|x|^{1-\alpha+\kappa})\\
&=C_\pm(\alpha,\kappa)|x|^{1-\alpha}[1+O(|x|^{\kappa})],\qquad x\to\pm 0,
\end{aligned}$$

and (6.12) holds.

*Case* $1\le\alpha<2$. In this case we use the Fourier inversion formula for distribution functions to get

$$\int_{|x|<\varepsilon}\rho(x)\,dx=\frac{2\zeta_\rho}{\pi}\int_0^{\infty}\frac{\sin(\varepsilon u)}{u}\operatorname{Re}[\widetilde{\psi}(u)]\,du.$$

The representation

$$\begin{aligned}
&\int_0^{\infty}\frac{\sin(\varepsilon u)}{u}\operatorname{Re}[\widetilde{\psi}(u)]\,du\\
&\quad=\int_0^{a}\frac{\sin(\varepsilon u)}{u}\operatorname{Re}[\widetilde{\psi}(u)]\,du+\int_a^{\infty}\frac{\sin(\varepsilon u)}{u}\operatorname{Re}[\widetilde{\psi}(u)]\,du\\
&\quad=: I_1+I_2
\end{aligned}$$

and the asymptotic relation

$$\begin{aligned}
I_2&=C_+(\alpha,\kappa)\int_a^{\infty}\frac{\sin(u\varepsilon)}{u}u^{\alpha-2}\,du+O(\varepsilon^{2-\alpha+\kappa})\\
&=C_+(\alpha,\kappa)\varepsilon^{2-\alpha}[1+O(\varepsilon^{\kappa})],\qquad \varepsilon\to+0,
\end{aligned}$$

lead now to (6.12) provided that $\kappa\le\alpha-1$.

**A.3. Proof of Theorem 6.4.** The representation

$$\widetilde{\alpha}_U-\alpha=\int_0^{\infty}w^U(u)(\widetilde{\mathcal{Y}}(u)-\mathcal{Y}(u))\,du+R_U$$

and Lemma 6.1 imply that

$$\begin{aligned}\mathrm{E}|\widetilde{\alpha}_U-\alpha|^2 &\le 3\mathrm{E}\biggl[\int_0^\infty w^U(u)\zeta_1(u)\Delta(u)\,du\biggr]^2 \\ &\quad + 3\mathrm{E}\biggl[\int_0^\infty w^U(u)\zeta_2(u)\Delta^2(u)\,du\biggr]^2 + 3|R_U|^2.\end{aligned} \tag{A.1}$$

Let us consider the first term in (A.1);

$$\begin{aligned}&\mathrm{E}\biggl[\int_0^\infty w^U(u)\zeta_1(u)\Delta(u)\,du\biggr]^2 \\ &\quad = \biggl[\int_0^\infty w^U(u)\zeta_1(u)\mathrm{E}[\Delta(u)]\,du\biggr]^2 + \mathrm{Var}\biggl[\int_0^\infty w^U(u)\zeta_1(u)\Delta(u)\,du\biggr].\end{aligned}$$

Since

$$\zeta_1(u) = 2^{-1}|\phi(u)|^{-2}\log^{-1}(|\phi(u)|) = e^{2\eta|u|^\alpha \operatorname{Re}\tau(u)}/(2\eta|u|^\alpha \operatorname{Re}\tau(u)),$$

we have

$$\begin{aligned}&\int_0^\infty w^U(u)\zeta_1(u)\mathrm{E}[\Delta(u)]\,du \\ &\quad = \int_0^1 w^1(u)\zeta_1(Uu)\mathrm{E}[\Delta(Uu)]\,du \\ &\quad = U^{-\alpha}\int_0^1 \frac{w^1(u)e^{2\eta U^\alpha u^\alpha \operatorname{Re}\tau(Uu)}}{2\eta u^\alpha \operatorname{Re}\tau(Uu)}\mathrm{E}[\Delta(Uu)]\,du.\end{aligned} \tag{A.2}$$

Due to the localization principle (Laplace method) and the identity

$$\mathrm{E}[\Delta(u)] = \mathrm{E}|\widetilde{\phi}(u)|^2 - |\phi(u)|^2 = \varepsilon(1-|\phi(u)|^2),$$

the integral in (A.2) is asymptotically (as $U\to\infty$) less than or equal to

$$A\varepsilon U^{-\alpha}\int_{1-\delta}^1 w^1(u)u^{-\alpha}e^{2\eta U^\alpha u^\alpha}\,du \lesssim \varepsilon U^{-\alpha}e^{2\eta U^\alpha}$$

with arbitrary small $\delta>0$ and some constant $A>0$. Similarly,

$$\begin{aligned}&\mathrm{Var}\biggl[\int_0^\infty w^U(u)\zeta_1(u)\Delta(u)\,du\biggr] \\ &\quad = \int_0^\infty\int_0^\infty w^U(u)w^U(v)\zeta_1(u)\zeta_1(v)\operatorname{Cov}(\Delta(u),\Delta(v))\,du\,dv \\ &\quad \lesssim \varepsilon U^{-2\alpha}e^{2\eta U^\alpha} + \varepsilon^2 U^{-4\alpha}e^{4\eta U^\alpha}, \qquad U\to\infty,\end{aligned}$$

where again localization principle and the identity,

$$\begin{aligned}
&\mathrm{Cov}(|\widetilde{\phi}(u)|^2, |\widetilde{\phi}(v)|^2)\\
&\quad = 2\varepsilon^3(\varepsilon^{-1}-1)(\varepsilon^{-1}-2)\\
&\qquad \times [\mathrm{Re}(\phi(u)\phi(v)\phi(-u-v))\\
&\qquad\quad + \mathrm{Re}(\phi(-u)\phi(v)\phi(u-v)) - 2|\phi(u)|^2|\phi(v)|^2]\\
&\qquad + \varepsilon^3(\varepsilon^{-1}-1)[|\phi(u+v)|^2 + |\phi(-u+v)|^2 - 2|\phi(u)|^2|\phi(v)|^2],
\end{aligned}$$

are used. Turn now to the second term in (A.1);

$$\begin{aligned}
&\mathrm{E}\left[\int_0^\infty w^U(u)\zeta_2(u)\Delta^2(u)\,du\right]^2\\
&\quad = \left[\int_0^\infty w^U(u)\zeta_2(u)\mathrm{E}[\Delta^2(u)]\,du\right]^2 + \mathrm{Var}\left[\int_0^\infty w^U(u)\zeta_2(u)\Delta^2(u)\,du\right].
\end{aligned}$$

Since

$$\zeta_2(u) \lesssim \frac{|\log|\phi(u)||}{|\phi(u)|^4}, \qquad u\to\infty,$$

and

$$\begin{aligned}
\mathrm{E}||\widetilde{\phi}(u)|^2 - |\phi(u)|^2|^2 &= \mathrm{E}||\widetilde{\phi}(u)|^2 - \mathrm{E}|\widetilde{\phi}(u)|^2 + \mathrm{E}|\widetilde{\phi}(u)|^2 - |\phi(u)|^2|^2\\
&\le 2\mathrm{E}||\widetilde{\phi}(u)|^2 - \mathrm{E}|\widetilde{\phi}(u)|^2|^2 + 2|\mathrm{E}|\widetilde{\phi}(u)|^2 - |\phi(u)|^2|^2\\
&\lesssim \varepsilon|\phi(u)|^2 + \varepsilon^2, \qquad u\to\infty,
\end{aligned}$$

we get an asymptotic estimate;

$$\int_0^\infty w^U(u)\zeta_2(u)\mathrm{E}[\Delta^2(u)]\,du \lesssim \varepsilon U^\alpha e^{2\eta U^\alpha} + \varepsilon^2 U^\alpha e^{4\eta U^\alpha}, \qquad U\to\infty.$$

Similarly, one can prove that

$$\mathrm{Var}\left[\int_0^\infty w^U(u)\zeta_2(u)\Delta^2(u)\,du\right] \lesssim \varepsilon^2 U^{2\alpha} e^{4\eta U^\alpha}, \qquad U\to\infty.$$

Finally, the third term in (A.1),

$$R_U = \int_0^\infty w^U(u)\log(\mathrm{Re}\,\tau(u))\,du,$$

can be can be bounded by

$$\begin{aligned}
|R_U| &= \left|\int_0^1 w^1(u)\log(\mathrm{Re}\,\tau(uU))\,du\right|\\
&\le U^{-1}\int_0^A |w^1(y/U)||\log(\mathrm{Re}\,\tau(y))|\,dy
\end{aligned}$$

$$+ U^{-\varkappa} \int_0^1 |y|^{-\varkappa} |w^1(y)| \, dy \lesssim U^{-\varkappa}, \qquad U \to \infty,$$

for $A > 0$ large enough. Combining all the previous estimates we get

$$\begin{aligned} \mathrm{E}|\widetilde{\alpha}_U - \alpha|^2 &\lesssim \varepsilon U^{-2\alpha} e^{2\eta U^{\alpha}} + \varepsilon^2 U^{2\alpha} e^{4\eta U^{\alpha}} + U^{-2\varkappa} \\ &\lesssim \varepsilon U^{-2\overline{\alpha}} e^{2\eta_+ U^{\overline{\alpha}}} + \varepsilon^2 U^{2\overline{\alpha}} e^{4\eta_+ U^{\overline{\alpha}}} + U^{-2\varkappa}, \qquad U \to \infty. \end{aligned} \tag{A.3}$$

Finally the choice

$$U = \left[ \frac{1}{2\eta_+} \log(\varepsilon^{-1} \log^{-\beta}(1/\varepsilon)) \right]^{1/\overline{\alpha}}$$

with $\beta = 1 + \varkappa/\overline{\alpha}$ leads to (6.13).

In the case of the calibration problem we have

$$\begin{aligned} |\widetilde{\phi}(u)|^2 = 1 - 2\,\mathrm{Re}\left[ u(u+\mathrm{i}) \sum_{j=1}^n \delta_j \widetilde{\mathcal{O}}(y_j) e^{\mathrm{i}uy_j} \right] \\ + u^2(1+u^2) \sum_{j,l=1}^n e^{\mathrm{i}u(y_l - y_j)} \delta_j \delta_l \widetilde{\mathcal{O}}(y_j) \widetilde{\mathcal{O}}(y_l) \end{aligned}$$

and

$$\begin{aligned} \mathrm{E}|\widetilde{\phi}(u)|^2 = 1 - 2\,\mathrm{Re}\left[ u(u+\mathrm{i}) \sum_{j=1}^n \delta_j \widetilde{O}(y_j) e^{\mathrm{i}uy_j} \right] \\ + u^2(1+u^2) \sum_{j \neq l}^n e^{\mathrm{i}u(y_l - y_j)} \delta_j \delta_l \widetilde{O}(y_j) \widetilde{O}(y_l) \\ + u^2(1+u^2) \sum_{j=1}^n \delta_j^2 \widetilde{\sigma}_j^2. \end{aligned}$$

As was mentioned in Section 3.2, function $\widetilde{O}(y) = e^{-y} O(y)$ is nonnegative, Lipschitz and satisfies the Cramér condition,

$$\int_{\mathbb{R}} O(y) e^{-y} \, dy < \infty,$$

provided that $\mathrm{E}[e^{2Y_T}] < \infty$. Under the condition $e^{-A} \le \|\delta\|^2$ we get

$$\left| \int_{\mathbb{R}} e^{\mathrm{i}uy} \widetilde{O}(y) \, dy - \sum_{j=1}^n e^{\mathrm{i}uy_j} \delta_j \widetilde{O}(y_j) \right| \lesssim \|\delta\|^2, \qquad \|\delta\|^2 \to 0,$$

as well as

$$\left|\left|\int_{\mathbb{R}} e^{\mathrm{i}uy}\widetilde{O}(y)\,dy\right|^2 - \sum_{j,l=1}^{n} e^{\mathrm{i}u(y_l-y_j)}\delta_j\delta_l\widetilde{O}(y_j)\widetilde{O}(y_l)\right| \lesssim \|\delta\|^2, \qquad \|\delta\|^2 \to 0.$$

Thus,

$$|\mathrm{E}|\widetilde{\phi}(u)|^2 - |\phi(u)|^2| \lesssim u^2(1+u^2)\sum_{j=1}^{n}\delta_j^2(1+\widetilde{\sigma}_j^2).$$

Further,

$$\begin{aligned} |\widetilde{\phi}(u)|^2 - \mathrm{E}|\widetilde{\phi}(u)|^2 = -2\,\mathrm{Re}\left[u(u+\mathrm{i})\sum_{j=1}^{n}\delta_j\widetilde{\sigma}_j\xi_j e^{\mathrm{i}uy_j}\right] \\ + 2u^2(1+u^2)\sum_{j<l} e^{\mathrm{i}u(y_l-y_j)}\delta_j\delta_l\widetilde{\sigma}_j\widetilde{\sigma}_l\xi_j\xi_l \\ + u^2(1+u^2)\sum_{j=1}^{n}\delta_j^2\widetilde{\sigma}_j^2(\xi_j^2-1) \end{aligned}$$

and

$$\mathrm{E}(|\widetilde{\phi}(u)|^2 - \mathrm{E}|\widetilde{\phi}(u)|^2)^2 \lesssim u^2(1+u^2)\sum_{j=1}^{n}\delta_j^2\widetilde{\sigma}_j^2 + u^4(1+u^2)^2\sum_{j=1}^{n}\delta_j^4\widetilde{\sigma}_j^4.$$

Using these inequalities, the first term in (A.1) can be estimated from above as

$$\begin{aligned} \mathrm{E}\left[\int_0^\infty w^U(u)\zeta_1(u)\Delta(u)\,du\right]^2 &\lesssim U^{8-2\alpha}e^{4\eta U^\alpha}\|\delta\|^4 + U^{4-2\alpha}e^{4\eta U^\alpha}\left[\sum_{j=1}^{n}\delta_j^2\widetilde{\sigma}_j^2\right]^2 \\ &\lesssim \varepsilon^2 U^{8-2\overline{\alpha}}e^{4\eta_+ U^{\overline{\alpha}}}, \end{aligned}$$

while the second one is asymptotically negligible if $\varepsilon^2 U^{8-2\alpha}e^{4\eta U^\alpha} \to 0$. Taking

$$U = \left[\frac{1}{2\eta_+}\log(\varepsilon^{-1}\log^{-\beta}(1/\varepsilon))\right]^{1/\overline{\alpha}}$$

with $\beta = (\varkappa+4)/\overline{\alpha} - 1$, we get (6.13).

**A.4. Proof of Theorem 6.5.** For any two probability measures $P$ and $Q$ define

$$\chi^2(P,Q) =: \begin{cases} \int \left(\frac{dP}{dQ} - 1\right)^2 dQ, & \text{if } P \ll Q, \\ +\infty, & \text{otherwise.} \end{cases}$$

The following proposition is the main tool for the proof of lower bounds in the estimation case and can be found in Butucea and Tsybakov (2004).

Proposition A.1. *Let $\mathcal{P}_\Theta := \{P_\theta : \theta \in \Theta\}$ be a family of models. Assume that there exist $\theta_1$ and $\theta_2$ in $\Theta$ with $|\theta_1 - \theta_2| > 2\delta_n > 0$ such that*

$$P_{\theta_1} \ll P_{\theta_2}, \qquad \chi^2(P_{\theta_1}^{\otimes n}, P_{\theta_2}^{\otimes n}) \leq \kappa^2 < 1,$$

*then*

$$\liminf_{n\to\infty} \inf_{\widehat{\theta}_n} \delta_n^{-2} \max\{\mathrm{E}_{\theta_1}|\widehat{\theta}_n - \theta_1|^2, \mathrm{E}_{\theta_2}|\widehat{\theta}_n - \theta_2|^2\} \geq (1-\kappa)^2(1-\sqrt{\kappa})^2,$$

*where the infimum is taken over all estimators $\widehat{\theta}_n$ (measurable function of observations) of the underlying parameter.*

Taking $\Theta = \mathcal{A}(\overline{\alpha}, \eta_-, \eta_+, \varkappa)$ and $\theta_i = (\alpha_i, \eta_i, \tau_i), i = 1, 2$, we get from Proposition A.1,

$$\sup_{(\alpha,\eta,\tau)\in\mathcal{A}(\overline{\alpha},\eta_-,\eta_+,\varkappa)} \mathrm{E}(|\alpha_\varepsilon - \alpha|^2) \geq \delta_n^{-2} \max\{\mathrm{E}_1(|\alpha_\varepsilon - \alpha_1|^2), \mathrm{E}_2(|\alpha_\varepsilon - \alpha_2|^2)\},$$

provided that $|\alpha_1 - \alpha_2| > 2\delta_n > 0$, and

$$\chi^2(P_{\theta_1}^{\otimes n}, P_{\theta_2}^{\otimes n}) \leq \kappa^2 < 1.$$

Turn now to the construction of models $\theta_1$ and $\theta_2$. Let us consider a symmetric stable model,

$$\psi(u) = \mathrm{i}\mu u + \vartheta(u), \qquad \vartheta(u) = -\eta_+|u|^\alpha, \qquad 0 < \alpha \leq 1, u \in \mathbb{R}.$$

For any $\delta$ satisfying $0 < \delta < \alpha$ and $M > 0$, define

$$\psi_\delta(u) = \mathrm{i}\mu u + \vartheta_\delta(u),$$

where

$$\vartheta_\delta(u) = -\eta_+|u|^\alpha \mathbf{1}_{\{|u|\leq M\}} - \frac{\eta_+ M^\delta}{(1 + cM^{-\varkappa})}|u|^{\alpha-\delta}(1 + c|u|^{-\varkappa})\mathbf{1}_{\{|u|>M\}}.$$

Then $\phi_\delta(u) = \exp(\mathrm{i}\mu u + \vartheta_\delta(u))$ is a characteristic function of some Lévy process and

$$\phi_\delta(u) = \phi(u), \qquad |u| \leq M,$$

where $\phi(u) = \exp(\mathrm{i}\mu u + \vartheta_\delta(u))$. Indeed, the function $\vartheta_\delta(u)$ is a continuous, nonpositive, symmetric function which is convex on $\mathbb{R}_+$ for large enough $M$ and small enough $c > 0$. According to a well-known Pólya criteria [see, e.g., Ushakov (1999)], the function $\exp(\xi\vartheta_\delta(u))$ is a cf. of some absolutely continuous distribution for any $\xi > 0$. In particular, for any natural $n$ the function $\exp(\vartheta_\delta(u)/n)$ is a cf. of some absolutely continuous distribution. Hence, $\exp(\vartheta_\delta(u))$ is a cf. of some infinitely divisible distribution. Define

$$\theta_1 = (\alpha, \eta_+, 1), \qquad \theta_2 = (\alpha - \delta, \eta_+, \tau_{\delta,M}) \tag{A.4}$$

and $\phi_{\theta_1}(u) = \phi(u), \phi_{\theta_2}(u) = \phi_\delta(u)$ with

$$\tau_{\delta,M}(u) := |u|^\delta \mathbf{1}_{\{|u|\le M\}} + \frac{M^\delta}{(1 + cM^{-\varkappa})}(1 + c|u|^{-\varkappa})\mathbf{1}_{\{|u|>M\}}.$$

If $M^\delta = 1 + cM^{-\varkappa}$, that is,

$$\delta = \log(1 + cM^{-\kappa})/\log M \asymp cM^{-\kappa}/\log M, \qquad M \to \infty, \tag{A.5}$$

then

$$|\tau_{\delta,M}(u) - 1| \lesssim |u|^{-\varkappa}, \qquad |u| \to \infty,$$

and hence $\theta_2 \in \Theta = \mathcal{A}(\overline{\alpha}, \eta_-, \eta_+, \varkappa)$. Furthermore, it holds

$$\chi^2(P_{\theta_1}^{\otimes n}, P_{\theta_2}^{\otimes n}) = n\chi^2(p_{\theta_1}, p_{\theta_2}) = n \int_{\mathbb{R}} \frac{|p_{\theta_1}(y) - p_{\theta_2}(y)|^2}{p_{\theta_1}(y)}\, dy,$$

where $p_{\theta_1}$ and $p_{\theta_2}$ are densities corresponding to cf. $\phi_{\theta_1}$ and $\phi_{\theta_2}$, respectively. Using the fact that the density of stable law $p_{\theta_1}(y)$ does not vanish on any compact set in $\mathbb{R}$ and fulfills

$$p_{\theta_1}(y) \gtrsim |y|^{-(\alpha+1)}, \qquad |y| \to \infty,$$

we derive

$$\begin{aligned} n\chi^2(p_{\theta_1}, p_{\theta_2}) &\le nC_1 \int_{|y|\le A} |p_{\theta_1}(y) - p_{\theta_2}(y)|^2\, dy \\ &\quad + nC_2 \int_{|y|>A} |y|^{\alpha+1} |p_{\theta_1}(y) - p_{\theta_2}(y)|^2\, dy \\ &= nC_1 I_1 + nC_2 I_2 \end{aligned}$$

for large enough $A > 0$ and some constants $C_1, C_2 > 0$. Using the fact that function $\phi_{\theta_1}(u) - \phi_{\theta_2}(u)$ is two times differentiable (it is zero for $|u| < M$) and Parseval's identity, we get

$$\begin{aligned} I_1 &\le \frac{1}{2\pi} \int_{\mathbb{R}} |\phi_{\theta_1}(u) - \phi_{\theta_2}(u)|^2\, du \\ &\le \frac{1}{2\pi} \int_{|u|>M} e^{-2\eta|u|^{\alpha-\delta}}\, du \lesssim M^{1-\alpha+\delta} e^{-2\eta M^{\alpha-\delta}}, \\ I_2 &\le \frac{1}{2\pi} \int_{|u|>M} |(\phi_{\theta_1}(u) - \phi_{\theta_2}(u))''|^2\, du \\ &\lesssim \int_{|u|>M} |u|^6 e^{-2\eta|u|^{\alpha-\delta}}\, du \lesssim M^{7-\alpha+\delta} e^{-2\eta M^{\alpha-\delta}}. \end{aligned}$$

The choice $M \asymp [\frac{1}{2\eta_+}\log(\varepsilon^{-1}\log^{-\beta}(1/\varepsilon))]^{1/(\alpha-\delta)}$ with $\varepsilon = 1/n$ and some $\beta \geq (7-(\alpha-\delta))/2(\alpha-\delta)$ yields

$$\varepsilon^{-1}\chi^2(p_{\theta_1}, p_{\theta_2}) < 1$$

for small enough $\varepsilon$. Combining this and (A.5), we arrive at (6.14).

For the proof of lower bounds in the case of calibration, one can employ the fact that the regression model,

$$\widetilde{\mathcal{O}}_T(y_i) = \widetilde{O}_T(y_i) + \widetilde{\sigma}(y_i)\xi_i, \qquad \delta_i = y_i - y_{i-1},$$
$$\mathrm{E}[\xi_i^2] = 1, \qquad i = 1, \ldots, n,$$

is equivalent to the Gaussian white noise model,

$$dZ(x) = \widetilde{O}(y)\,dy + \varepsilon^{1/2}\,dW(y)$$

with the noise level asymptotics $\varepsilon \to 0$, a two-sided Brownian motion $W$. Here the noise level $\varepsilon$ corresponds to the statistical regression error $\sum_{j=1}^n \delta_j^2 \widetilde{\sigma}_j^2$. Furthermore, instead of $\chi^2$ distance we use the Kullback–Leibler divergence,

$$\mathrm{KL}(\mathcal{T}_{\theta_1}, \mathcal{T}_{\theta_2}) = \frac{1}{2}\int_{\mathbb{R}} |(\widetilde{O}_{\theta_1} - \widetilde{O}_{\theta_2})(y)|^2 \varepsilon^{-1}\,dy,$$

between two models $\mathcal{T}_{\theta_1}$ and $\mathcal{T}_{\theta_2}$ corresponding to two Lévy processes with characteristics $\theta_1$ and $\theta_2$, respectively [see (A.4)]. Simple calculations lead to the estimate,

$$\mathrm{KL}(\mathcal{T}_{\theta_1}, \mathcal{T}_{\theta_2}) \lesssim \varepsilon^{-1} M^{\gamma} e^{-2\eta_+ M^{\alpha-\delta}}$$

with some $\gamma > 0$. Hence, for small enough $\varepsilon > 0$ it holds

$$\mathrm{KL}(\mathcal{T}_{\theta_1}, \mathcal{T}_{\theta_2}) < 1$$

provided that $M \asymp [\frac{1}{2\eta_+}\log(\varepsilon^{-1}\log^{-\beta}(1/\varepsilon))]^{1/(\alpha-\delta)}$ with $\beta \geq \gamma/2(\alpha-\delta)$. The Assouad lemma [see, e.g., Tsybakov (2008)] together with (A.5) implies (6.14).

**A.5. Proof of Proposition 6.6.** It holds for any fixed $U$,

$$\begin{aligned}
\widetilde{\alpha}_U - \alpha &= \int_0^\infty w^U(u)(\widetilde{\mathcal{Y}}(u) - \mathcal{Y}(u))\,du \\
&= \int_0^\infty w^U(u)\zeta_1(u)\Delta(u)\,du \\
&\quad + \int_0^\infty w^U(u)Q(u)\,du + R_U,
\end{aligned}$$

where $Q$ is defined in (6.3). As shown in Lemma A.2, the process $\varepsilon^{-1/2}\Delta(u)$ converges weakly to a Gaussian process $Z(u)$ with $\mathrm{E}[Z(u)] = 0$ and $\mathrm{Cov}(Z(u),$ $Z(v)) = S(u,v)$. Moreover, $\varepsilon^{-1/2}Q(u) \to 0$, almost surely. The representation for $\delta(u)$ in Lemma A.2 and CLT for $U$-statistics implies that if for some sequence $U(\varepsilon)$, it holds $\varsigma^{-1}(\varepsilon)R_{U(\varepsilon)} \to 0$ with

$$\varsigma^2(\varepsilon) = \varepsilon \int_0^\infty w^{U(\varepsilon)}(u) w^{U(\varepsilon)}(v) \zeta_1(u)\zeta_1(v) S(u,v)\, du\, dv,$$

then $\varsigma^{-1}(\varepsilon)(\widetilde{\alpha}_{U(\varepsilon)} - \alpha) \to \mathcal{N}(0,1)$.

**A.6. Proof of Theorem 6.7.** We give only the sketch of the proof. Let $\omega_-$ and $\omega_+$ be two truncation levels satisfying $0 < \omega_-(u) < \rho_\xi(u) < \omega_+(u) < 1$ and $0 < \omega_- < \rho_\xi(u)(1 - \log(\rho_\xi(u))/(1 + \log(\rho_\xi(u))))$. First, similar to the proof of Proposition 6.1, one can show that

$$|\widetilde{\mathcal{Y}}_\xi(u) - \mathcal{Y}(u) - \zeta_{1,\xi}(u)(T_{0,\omega_+}[\widetilde{\rho}_\xi](u) - \rho_\xi(u))| \le \zeta_{2,\xi}(u)(T_{0,\omega_+}[\widetilde{\rho}_\xi](u) - \rho_\xi(u))^2,$$

where

$$\zeta_{1,\xi}(u) = -\rho_\xi^{-1}(u) \log^{-1}(\rho_\xi(u))$$

and

$$\zeta_2(u) = \max_{\theta \in \{\omega_-(u), \omega_+(u)\}} \left[ \frac{1 + |\log(\theta)|}{\theta^2 \log^2(\theta)} \right].$$

Furthermore, we have on the set $\{\widetilde{\rho}_\xi(u) \le \omega_+(u)\}$

$$\begin{aligned}
&|\rho_\xi(u) - T_{0,\omega_+}[\widetilde{\rho}_\xi](u)| \\
&\qquad \le \omega_+(u) \frac{||\phi_a(\xi u)|^2 - |\widetilde{\phi}(\xi u)|^2|}{|\phi_a(\xi u)|^2} + \frac{||\phi_a(u)|^{2\xi^2} - |\widetilde{\phi}(u)|^{2\xi^2}|}{|\phi_a(\xi u)|^2},
\end{aligned}$$

and on the set $\{\widetilde{\rho}_\xi(u) > \omega_+(u)\}$ it holds

$$|\rho_\xi(u) - T_{0,\omega_+}[\widetilde{\rho}_\xi](u)| \le 2\omega_+(u).$$

Hence

$$\begin{aligned}
&\mathrm{E}|\rho_\xi(u) - T_{0,\omega_+}[\widetilde{\rho}_\xi](u)|^2 \\
&\qquad \le 2|\phi_a(\xi u)|^{-4}[\mathrm{E}||\phi_a(\xi u)|^2 - |\widetilde{\phi}(\xi u)|^2|^2 \\
&\qquad\qquad + \mathrm{E}||\phi_a(u)|^{2\xi^2} - |\widetilde{\phi}(u)|^{2\xi^2}|^2] + 4\omega_+^2(u)\,\mathrm{P}(\widetilde{\rho}_\xi(u) > \omega_+(u)).
\end{aligned}$$

Without loss of generality one can assume that there exists $U_0 > 0$ such that $\rho_\xi(u)/\omega_+(u) < 1/2$ for $u > U_0$. Then it holds for $u > U_0$

$$\begin{aligned}
&\mathrm{P}(\widetilde{\rho}_\xi(u) > \omega_+(u)) \\
&\qquad \le \mathrm{P}(||\phi_a(u)|^{2\xi^2} - |\widetilde{\phi}(u)|^{2\xi^2}| > \omega_+(u)|\phi(u\xi)|^2/4) \\
&\qquad\quad + \mathrm{P}(||\phi_a(u\xi)|^2 - |\widetilde{\phi}(u\xi)|^2| > \omega_+(u)|\phi(u\xi)|^2/4) \\
&\qquad \le 16|\phi_a(\xi u)|^{-4}[\mathrm{E}||\phi_a(\xi u)|^2 - |\widetilde{\phi}(\xi u)|^2|^2 + \mathrm{E}||\phi_a(u)|^{2\xi^2} - |\widetilde{\phi}(u)|^{2\xi^2}|^2].
\end{aligned}$$

In the case of the estimation under $\mathbb{P}$, for instance, we have

$$\mathrm{E}||\phi_a(\xi u)|^2 - |\widetilde{\phi}(\xi u)|^2|^2 \lesssim \varepsilon, \qquad \mathrm{E}||\phi_a(u)|^{2\xi^2} - |\widetilde{\phi}(u)|^{2\xi^2}|^2 \lesssim \varepsilon, \qquad \varepsilon \to 0,$$

and hence

$$\mathrm{E}|\rho_\xi(u) - T_{0,\omega_+}[\widetilde{\rho}_\xi](u)|^2 \lesssim \varepsilon |\phi_a(\xi u)|^{-4}, \qquad \varepsilon \to 0.$$

Now one can follow the proof of Theorem 6.4 and use the fact that

$$\zeta_{1,\xi}(u) \asymp c_\xi^{-1}(\alpha)|u|^{-\alpha}\tau_\xi^{-1}(u)\exp(c_\xi(\alpha)|u|^\alpha \tau_\xi(u)), \qquad u \to \infty.$$

**A.7. Proof of Theorem 6.8.** Instead of Lévy models $\theta_1$ and $\theta_2$, one considers models $\theta_{1,a}$ and $\theta_{2,a}$ with characteristic exponents $\psi_a(u) = \mathrm{i}\mu u - \overline{a}^2 u^2/2 + \vartheta(u)$ and $\psi_{a,\delta}(u) = \mathrm{i}\mu u - \overline{a}^2 u^2/2 + \vartheta_\delta(u)$, respectively. The rest of the proof is almost identical to the proof of Theorem 6.5.

**A.8. Auxiliary results.** The following lemma is a basic tool to investigate the asymptotic behavior of the estimate $\widetilde{\alpha}$ under the historical measure $\mathbb{P}$.

Lemma A.2. *The process $\varepsilon^{-1/2}\Delta(u)$ with $\Delta(u) = |\widetilde{\phi}(u)|^2 - |\phi(u)|^2$ weakly converges to a Gaussian process $Z(u)$ with $\mathrm{E}[Z(u)] = 0$ and $\mathrm{Cov}(Z(u), Z(v)) = S(u,v)$ where*

$$\begin{aligned} S(u,v) &:= \operatorname{Re}\phi(u-v) + \operatorname{Im}\phi(u+v) \\ &\quad - (\operatorname{Re}\phi(u) + \operatorname{Im}\phi(u))(\operatorname{Re}\phi(v) + \operatorname{Im}\phi(v)). \end{aligned}$$

Proof. We have

$$|\widetilde{\phi}(u)|^2 = [\operatorname{Re}\widetilde{\phi}(u)]^2 + [\operatorname{Im}\widetilde{\phi}(u)]^2 = \frac{1}{n^2}\sum_{j=1}^n \sum_{k=1}^n \cos(u(X_j - X_k)).$$

Put

$$H_n(u) = \binom{n}{2}^{-1} \sum_c \cos(u(X_j - X_k)) = \frac{2}{n(n-1)} \sum_c \cos(u(X_j - X_k)),$$

where summation $c$ is over all $\binom{n}{2}$ combinations of 2 integers chosen from $(1,\dots,n)$. Then

$$\varepsilon^{-1/2}(|\widetilde{\phi}(u)|^2 - |\phi(u)|^2) = \varepsilon^{1/2} + \varepsilon^{-1/2}(1-\varepsilon)(H_n - |\phi(u)|^2) - \varepsilon^{1/2}|\phi(u)|^2.$$

The first and third terms on the right-hand side converge to 0. Consider the middle term. Since $H_n(u)$ is an $U$-statistic (for each $u$), $\varepsilon^{-1/2}(H_n - |\phi(u)|^2)$ weakly converges to a Gaussian process with zero mean and covariance

$$\mathrm{Cov}[\mathrm{E}_{X_2}\cos(u(X_1 - X_2)), \mathrm{E}_{X_2}\cos(v(X_1 - X_2))]$$

(where $\mathrm{E}_X Y$ denotes the conditional expectation of $Y$ given $X$). Let us compute this covariance. For any $u,v\in\mathbb{R}$ it holds

$$\begin{aligned}
&\operatorname{Cov}(\mathrm{E}_{X_2}[\cos(u(X_1-X_2))],\mathrm{E}_{X_2}[\cos(v(X_1-X_2))])\\
&\quad=\mathrm{E}[(\cos(uX_2)-\operatorname{Re}\phi(u))\operatorname{Re}\phi(u)+(\sin(uX_2)-\operatorname{Im}\phi(u))\operatorname{Im}\phi(u)]\\
&\qquad\times[(\cos(vX_2)-\operatorname{Re}\phi(v))\operatorname{Re}\phi(v)+(\sin(vX_2)-\operatorname{Im}\phi(v))\operatorname{Im}\phi(v)],
\end{aligned}$$

where

$$\begin{aligned}
&\mathrm{E}(\cos(uX_2)-\operatorname{Re}\phi(u))(\cos(vX_2)-\operatorname{Re}\phi(v))\\
&\quad=\frac{\operatorname{Re}\phi(u+v)+\operatorname{Re}\phi(u-v)}{2}-\operatorname{Re}\phi(u)\operatorname{Re}\phi(v),\\
&\mathrm{E}(\sin(uX_2)-\operatorname{Im}\phi(u))(\sin(vX_2)-\operatorname{Im}\phi(v))\\
&\quad=\frac{\operatorname{Re}\phi(u-v)-\operatorname{Re}\phi(u+v)}{2}-\operatorname{Im}\phi(u)\operatorname{Im}\phi(v)
\end{aligned}$$

and

$$\begin{aligned}
&\mathrm{E}(\cos(uX_2)-\operatorname{Re}\phi(u))(\sin(vX_2)-\operatorname{Im}\phi(v))\\
&\quad=\frac{\operatorname{Im}\phi(v-u)+\operatorname{Im}\phi(u+v)}{2}-\operatorname{Re}\phi(u)\operatorname{Im}\phi(v).
\end{aligned}$$ □

Weierstrass-Institute
Mohrenstr. 39
10117 Berlin
Germany
E-mail: belomest@wias-berlin.de